\documentclass[10pt]{amsproc}
\usepackage{amsmath,amsthm,amsfonts,amscd,float}
\newcommand{\Gal}{{\rm Gal}}
\newcommand{\XT}{\tilde{X}}
\newcommand{\ST}{\tilde{S}}
\newcommand{\CC}{\mathbb{C}}
\newcommand{\PP}{\mathbb{P}}
\newcommand{\FF}{\mathbb{F}}
\newcommand{\OO}{\mathcal{O}}
\newcommand{\ZZ}{\mathbb{Z}}
\newcommand{\QQ}{\mathbb{Q}}

\newcommand{\Frob}{{\rm Frob}}

\newcommand{\image}{{\rm im}\thinspace}
\newcommand{\Trace}{{\rm tr}\thinspace}
\newcommand{\trace}{{\rm tr}\thinspace}

\newcommand{\Ind}{{\rm Ind}\thinspace}
\newcommand{\Mod}{\negthinspace{\rm mod}\thinspace}
\newcommand{\Fpbar}{\overline{\mathbb{F}}_p}
\newcommand{\Qpbar}{\overline{\mathbb{Q}}_p}
\newcommand{\Qbar}{\overline{\mathbb{Q}}}
\newcommand{\kbar}{\overline{k}}

\theoremstyle{plain}
\newtheorem{thm}{Theorem}[section]

\newtheorem{lem}[thm]{Lemma}

\newtheorem{prop}[thm]{Proposition}

\theoremstyle{definition}
\newtheorem{defn}{Definition}[section]
\newtheorem{conj}{Conjecture}[section]

\newtheorem*{case1}{Case 1}
\newtheorem*{case2}{Case 2}
\newtheorem*{step1}{Step 1}
\newtheorem*{step2}{Step 2}
\newtheorem*{step3}{Step 3}
\newtheorem*{question}{Question}

\theoremstyle{definition}
\newtheorem*{rem}{Remark}

\floatstyle{plain}
\newfloat{Program}{thp}{lop}[section]
\newfloat{Table}{thp}{lop}[section]

\begin{document}

\title{A Modular Non-Rigid Calabi-Yau Threefold}
\author{Edward Lee}
\address{Department of Mathematics, Harvard University, Cambridge, MA 02138}
\email{edward.lee@post.harvard.edu}


 
\begin{abstract}
We construct an algebraic variety by resolving singularities of a 
quintic Calabi-Yau threefold.  The middle cohomology of the threefold
is shown to contain a piece coming from a pair of elliptic surfaces.
The resulting quotient is a two-dimensional Galois representation.  By
using the Lefschetz fixed-point theorem in \'etale cohomology and
counting points on the variety over finite fields, this Galois
representation is shown to be modular.
\end{abstract}
\maketitle


\section{Introduction}

In this paper we investigate the geometry and arithmetic of a
Calabi-Yau threefold $X \subset \PP^4 \times \PP^4$ given as a
complete intersection of five hypersurfaces of bidegree $(1,1)$.
$X$ is a partial desingularization common to a pair of quintic
threefolds $F$ and $G$ in $\PP^4$; we will in fact be interested in a
big resolution $\XT$ of $X$.  So $\XT$ is not actually a
Calabi-Yau threefold, but it is birational to one.

A conjecture of Fontaine and Mazur \cite{bib:FM} predicts that
two-dimensional $l$-adic Galois representations coming from geometry
should be modular.  More precisely, the statement is
that a continuous irreducible two-dimensional $l$-adic representation
of the absolute Galois group $G_{\QQ}$ that is isomorphic to a Tate
twist of a subquotient of an \'etale cohomology group of a variety
$X/\QQ$ should be modular.  This is a higher-dimensional
generalization of the Taniyama-Shimura conjecture on the modularity of
elliptic curves over $\QQ$, as proved by Wiles, Taylor et al. \cite{bib:TW},
\cite{bib:BCDT}.

A rigid Calabi-Yau threefold $X$ defined over $\QQ$ has two-dimensional
middle cohomology, and is thus expected to be modular.  We expect the
$L$-series of the Galois action on the $l$-adic cohomology to be, up
to factors associated to the primes of bad reduction of $X$, the
Mellin transform of a weight 4 modular form.  Recently Dieulefait and
Manoharmayum \cite{bib:DM} proved that rigid Calabi-Yau threefolds
that have good reduction at 3 and 7, or at 5 and another suitable
prime, are modular.  A handful of explicit examples of modular
Calabi-Yau threefolds are known.  Some examples are given in
\cite{bib:Schoen}, \cite{bib:WvG}, \cite{bib:HV}, \cite{bib:Y1},
\cite{bib:Y2}.

Some nonrigid Calabi-Yau threfolds have been shown to be
modular by Hulek-Verrill \cite{bib:HV} and Sch\"utt \cite{bib:Schutt}.
Here the middle cohomology group of $X$ has dimension greater than 2,
so one must figure out how to extract a 2-dimensional piece on which
$G_\QQ$ acts.  

Hulek-Verrill found threefolds in a toric variety with
$h^3 = 4, 6$ and $10$; in each case they showed that the semisimplification of $H^3$ was a
direct sum $\oplus_i W_i \oplus V$, where the $W_i$ came from
elliptic surfaces defined over $\QQ$ and $V$ was the
remaining quotient.  The $W_i$ were thus isomorphic to
cohomology groups of elliptic curves twisted by $(-1)$, and $V$ was
shown to correspond to a modular form of weight 4.

Sch\"utt constructed some fiber products of rational elliptic surfaces; he
showed that each of their middle cohomology groups also broke up into
a sum of two-dimensional pieces coming from elliptic curves and a leftover
two-dimensional piece corresponding to a modular form of weight 4.

Our threefold $\XT$ is constructed by resolving singularities of
a pair of Horrocks-Mumford quintic threefolds $F$ and $G$.  The Horrocks-Mumford bundle
$HM$ is a stable indecomposable rank 2 vector bundle over $\PP^4$, and
zero sets of its sections are abelian surfaces.  We take zero sets of
a pair of sections of the line bundle $\wedge^2(HM)
\cong \OO(5)$ as our quintic threefolds $F$ and $G$; they are thus
pencils of abelian surfaces.  

First we take a common partial resolution $X$ of $F$ and $G$ as
mentioned at the outset; it is a Calabi-Yau threefold given as a
complete intersection in $\PP^4 \times \PP^4$.  We then blow up the
singularities of $X$ to obtain $\XT$.  Unfortunately, there exists no model for a small resolution
of $X$ over $\QQ$.

By studying the geometry of $\XT$ and
by exploiting the Weil conjectures, we are able to show that $h^3(\XT)
= 6$.  We show that the semisimplification of the Galois
representation $H^3(\XT)$ is a direct sum of a two-dimensional piece
$V$ and a four-dimensional piece $W$.  The four-dimensional piece $W$
arises from the cohomology of a pair of elliptic surfaces $E_1$ and
$E_2$ that are complex conjugates of each other.  Thus their union is
defined over $\QQ$, and as a Galois representation $W$ is induced from
a representation of the subgroup $G_{\QQ(i)}$.  We then show that the
two-dimensional piece $V$ is modular; by using a theorem of
Faltings-Serre-Livn\'e \cite{bib:Livne}, we are able to prove this by
studying the reduction of $\XT$ modulo a finite set of primes.  In
practice this amounts to counting the points on $\XT$ over $\FF_p$, a
task which can easily be done by computer.

This paper is organized as follows:  In section 2 we review the
construction and key properties of the Horrocks-Mumford vector bundle.
In section 3 we construct the Horrocks-Mumford quintic threefolds $F$
and $G$ as pencils of abelian surfaces.  In section 4 we construct the
common resolution $\XT$ of the threefolds $F$ and $G$ and study its
geometry.  In section 5 we find the elliptic surfaces $E_1$ and $E_2$
in $\XT$, count points and apply Livn\'e's method to show that $\XT$ is modular.

This paper is based upon doctoral research conducted at Harvard
University.  The author would like to thank Shing-Tung Yau for his
guidance and encouragement, Richard Taylor and Klaus Hulek for
valuable discussions, and Jan Stienstra, Matthias Sch\"utt and the referee for some helpful comments.


\section{The Horrocks-Mumford vector bundle}

\subsection{Construction of the bundle}

The Horrocks-Mumford vector bundle $HM$ is a stable, indecomposable
rank 2 bundle over the complex projective space $\PP^4$.  It is
essentially the only known bundle satisfying these properties; all
other such bundles that are currently known are derived from $HM$ by
twisting by powers of the sheaf $\OO(1)$ or by taking pullbacks
to branched covers of $\PP^4$.  It was first discovered by Horrocks
and Mumford in \cite{bib:HM}, and has been further studied by many
other authors (see for example \cite{bib:BHM}, \cite{bib:BM},
\cite{bib:HL}, \cite{bib:Schoen}).  In this section we will describe the
construction of $HM$ and explain some of its properties that we will
use later.

The following exposition of the Horrocks-Mumford bundle is taken from
\cite{bib:Hulek}.

A {\em monad} is a three-term complex

\begin{equation*}
\begin{CD}
A @>p>> B @>q>> C
\end{CD}
\end{equation*}

\noindent of vector bundles where $p$ is injective and $q$ is surjective.  The
cohomology of the monad

\[
E = \ker q / \image p
\]

\noindent is also a vector bundle.  To construct the Horrocks-Mumford bundle
using a monad we fix a vector space

\[
V \cong \CC^5.
\]

Denote its standard basis by $e_i, i \in \ZZ/5$.  On the projective
space $\PP^4 = \PP(V)$, we have the Koszul complex

\begin{equation*}
\begin{split}
0 \longrightarrow \OO \stackrel{\wedge s}\longrightarrow V \otimes \OO(1)
&\stackrel{\wedge s}\longrightarrow {\wedge}^2 V \otimes \OO(2) \\
&\stackrel{\wedge s}\longrightarrow {\wedge}^3
V \otimes\OO(3) \stackrel{\wedge s}\longrightarrow {\wedge}^4 V
\otimes \OO(4) \longrightarrow \OO(5) \longrightarrow 0. \\
\end{split}
\end{equation*}

Now the quotient $\OO(1) \otimes V / \OO$ is isomorphic to the tangent
sheaf $T$, and in the Koszul complex the sheaf of cycles
$\image (\OO(i) \otimes \wedge^i V) \subset \OO(i+1) \otimes
\wedge^{i+1}V$ is isomorphic to $\wedge^i T$.  Thus from the map

\[
\OO(2) \otimes \wedge^2 V \longrightarrow \OO(3) \otimes \wedge^3 V
\]

\noindent we obtain the sequence of maps

\[
\OO(2) \otimes \wedge^2 V \stackrel{p_0}\longrightarrow \wedge^2 T \stackrel{q_0}\longrightarrow
\OO(3) \otimes \wedge^3 V
\]

\noindent where the first map is surjective and the second is
injective.

Horrocks and Mumford defined the following maps

\begin{equation*}
\begin{aligned}
f^+ &: V \longrightarrow \wedge^2 V, & f^+ (\sum v_i e_i) &= \sum v_i e_{i+2}
\wedge e_{i+3}, \\
f^- &: V \longrightarrow \wedge^2 V, & f^- (\sum v_i e_i) &= \sum v_i e_{i+1}
\wedge e_{i+4}. \\
\end{aligned}
\end{equation*}

Using these maps one can define

\begin{equation*}
\begin{aligned}
p&: V \otimes \OO(2) \stackrel{(f^+, f^-)(2)}\longrightarrow 2 \wedge^2 V
\otimes \OO(2) \stackrel{2 p_0}\longrightarrow 2 \wedge^2 T \\
q&: 2 \wedge^2 T \stackrel{2 q_0}\longrightarrow 2 \wedge^3 V \otimes \OO(3)
\stackrel{(-f^{-*}, f^{+*})(3)}\longrightarrow V^* \otimes \OO(3).
\end{aligned}
\end{equation*}

One easily checks that $q \circ p = 0$.  Hence we obtain a monad

\begin{equation*}
V \otimes \OO(2) \stackrel{p}\longrightarrow 2 \wedge^2 T
\stackrel{q}\longrightarrow V^* \otimes
\OO(3).
\end{equation*}

Its cohomology

\begin{equation*}
HM = \ker q / \image p
\end{equation*}

\noindent is the Horrocks-Mumford bundle.  It is a rank 2 bundle, and its total
Chern class $c(HM)$ equals $c(\wedge^2 T)^2 c(V^* \otimes \OO(3))^{-1}
c(V \otimes \OO(2))^{-1}$.  Using the splitting principle, one
computes this class to be $1 + 5H + 10 H^2$.  Therefore, zero sets of
sections of $HM$ are surfaces of degree 10; Horrocks and Mumford
showed that the generic zero set is a smooth abelian surface.

\subsection{Symmetries of $HM$ and invariant quintics}

The study of $HM$ has been greatly expedited by the fact that it
admits a large group of discrete symmetries.  Consider the Heisenberg
group of rank 5, which we denote by $H_5$.  We present it as a
subgroup of $GL_5(\CC)$ generated by the matrices

\[
\sigma = \begin{pmatrix} & 1 & & & \\ & & 1 & & \\ & & & 1 & \\ & & &
  & 1 \\ 1 & & & & \end{pmatrix}, \tau = \begin{pmatrix} 1 & & & & \\
  & \epsilon & & & \\ & & \epsilon^2 & & \\ & & & \epsilon^3 & \\ & &
  & & \epsilon^4 \end{pmatrix},
\]

\noindent where $\epsilon = e^{\frac{2 \pi i}{5}}$ is a primitive fifth root of
unity.  $H_5$ is a central extension

\[
1 \rightarrow \mu_5 \rightarrow H_5 \rightarrow \ZZ/5 \times \ZZ/5
\rightarrow 1
\]

\noindent where $\sigma$ is sent to $(1,0)$ and $\tau$ to $(0,1)$.
Here $\mu_5$ is the multiplicative group of fifth roots of unity.

In fact, the normalizer $N_5$ of $H_5$ in $SL_5(\CC)$ preserves $HM$.
$N_5$ is a semidirect product of $H_5$ with the binary icosahedral
group $SL(2,\ZZ_5)$.  We will need the following elements of $N_5$:

\[
\iota = \begin{pmatrix} 1 & & & & \\ & & & & 1 \\ & & & 1 & \\ & & 1 &
  & \\ & 1 & & & \end{pmatrix}, \mu = \begin{pmatrix} 1 & & & & \\ &
  & 1 & & \\ & & & & 1 \\ & 1 & & & \\ & & & 1 & \end{pmatrix},
  \nu = \begin{pmatrix} 1 & & & & \\ & \epsilon & & & \\ & &
  \epsilon^4 & & \\ & & & \epsilon^4 & \\ & & & & \epsilon
  \end{pmatrix}.
\]

These matrices act on $N_5/\mu_5 \simeq \ZZ/5 \times \ZZ/5$ by
conjugation; Horrocks and Mumford showed that the action is
unimodular.  Their images in $SL_2(\ZZ_5)$ are

\[
\overline{\iota} = \begin{pmatrix} -1 & \\ & -1 \end{pmatrix},
\overline{\mu} = \begin{pmatrix} 2 & \\ & 3 \end{pmatrix},
\overline{\nu} = \begin{pmatrix} 1 & 2 \\ & 1 \end{pmatrix}.
\]

En route to determining the sections of $HM$, Horrocks and Mumford
determined the $N_5/H_5$-module $\Gamma_{H_5}(\OO(5))$ of
$H$-invariants of $\Gamma(\OO(5))$, i.e. Heisenberg-invariant quintics
in $\PP^4$.  It is six-dimensional, spanned by the polynomials

\[
\sum x_i^5, \sum x_i^3 x_{i+1} x_{i+4}, \sum x_i x_{i+1}^2 x_{i+4}^2,
\]
\[
\sum x_i^3 x_{i+2} x_{i+3}, \sum x_i x_{i+2}^2 x_{i+3}^2, x_0 x_1 x_2
x_3 x_4
\]

\noindent where the sums are taken over powers of $\sigma$.  The base locus of
this space of quintics is the set of 25 lines $L_{ij}$,
where

\[
L_{00} = \{ x \in \PP^4 : x_0 = x_1 + x_4 = x_2 + x_3 = 0 \},
\]

\[
L_{ij} = \sigma^i \tau^j L_{00}.
\]

Since $c(HM) = 1 + 5H + 10H^2$, $c_1(\wedge^2(HM)) = 5H$ and thus
$\wedge^2(HM) \cong \OO(5)$.  Hence if $s_1$ and $s_2$ are sections of $HM$, the zero set
of the section $s_1 \wedge s_2$ of $\wedge^2{HM}$ is a (singular) quintic
Calabi-Yau threefold that has the structure of a pencil of abelian
surfaces.

\begin{prop}  For generic sections $s_1$ and $s_2$ of $HM$, the
singularities of the resulting threefold are the 100 nodes coming from
the intersection of $Z(s_1)$ and $Z(s_2)$.
\end{prop}

\begin{proof}  For a generic choice of $s_1$ and $s_2$,
$Z(s_1)$ and $Z(s_2)$ intersect transversely in 100
points; these points form the base locus of the pencil.

Let $p$ be a point at which $Z(s_1)$ and $Z(s_2)$ intersect
transversely.  Choose a trivialization of $HM$ near $p$, and put $s_1
= (s_{11}, s_{12})$ and $s_2 = (s_{21}, s_{22})$ relative to this
trivialization.  Since $Z(s_1)$ and $Z(s_2)$ intersect transversely, we may then use $s_{11}, s_{12}, s_{21}$ and $s_{22}$
as local coordinates on $\PP^4$ centered at $p$.  The local equation
for the threefold is then

\[
s_{11} s_{22} - s_{12} s_{21} = 0.
\]

\noindent Hence $p$ is a node.
\end{proof}

Nodes on threefolds result from the vanishing of an $S^3$ cycle on a a smooth 
family of threefolds.  One expects that degenerating the $S^3$ cycles and then 
resolving the singularities will cause the Betti number $h^3$ to drop.  We will 
be interested in birationally equivalent Calabi-Yau threefolds with low Betti number.  
Taking one-parameter families of abelian surfaces in $\PP^4$ gives us a quick way of
manufacturing nodal Calabi-Yau threefolds, whose singularities can
then be resolved.  In \cite{bib:Schoen}, Schoen studied the Fermat
quintic $Q$ defined by the equation 
$$x_0^5 + x_1^5 + x_2^5 + x_3^5 + x_4^5 - 5 x_0 x_1 x_2 x_3 x_4
= 0.$$  Schoen showed that it was a Horrocks-Mumford quintic with 125 nodes
instead of the usual 100, and he proved that the blowup $\tilde{Q}$ of
$Q$ was rigid and modular.  Other nodal Calabi-Yau threefolds whose
resolutions are modular were studied in \cite{bib:WvG}.

\begin{rem} Instead of manufacturing nodal Calabi-Yaus in $\PP^4$, one can also
use as the ambient space other Fano fourfolds such as $\PP^3 \times
\PP^1$; we then want to consider a rank 2 bundle whose determinant
bundle is anticanonical.  As above, we can then consider surfaces cut
out by sections of the bundle and take pencils of these surfaces to
obtain other nodal Calabi-Yau threefolds.  In \cite{bib:Lange} and
\cite{bib:Lange2}, Lange has proven the
existence of abelian surfaces in $\PP^1 \times \PP^3$ and by the Serre
construction found the rank 2 bundle $V$ whose zero sections yield these
surfaces.  
\end{rem}


\section{Abelian surfaces in $\PP^4$}

\subsection{Sections of $HM$}

In the previous section we mentioned that an abelian surface in $\PP^4$ is
projectively equivalent to $Z(s)$ for some section $s$ of $HM$.  Since
$h^0(HM)$ is 4, $\PP^3$ is a parameter space of (possibly degenerate)
abelian surfaces in $\PP^4$.

Any vector bundle over $\PP^1$ splits into a direct sum of line
bundles; for most lines in $\PP^4$, the restriction of $HM$ is
isomorphic to $\OO(2) \oplus \OO(3)$.  Lines $L$ such that $HM|_L$ is
isomorphic to $\OO(2-a) \oplus \OO(3+a)$ are called \textit{jumping lines} of
order $a$.  It is well known that the 25 lines $L_{ij}$ are jumping
lines of order $3$; the restriction of $HM$ to these lines is isomorphic to
$\OO(-1) \oplus \OO(6)$.  Barth, Hulek and Moore (\cite{bib:BHM}) proved that the
restriction map $\Gamma_{HM} \longrightarrow \Gamma_{HM|_{L_{00}}}$ is
injective, and they were able to determine the sections of
$\Gamma_{HM|_{L_{00}}}$:

\begin{prop}
Let $\lambda$ and $\mu$ be the restrictions of the coordinates $x_1$
and $x_2$ to $L_{00}$.  Then the image of $\Gamma_{HM}$ in
$\Gamma_{HM|_{L_{00}}}$ is spanned by the sections $t_0 = \lambda^6 +
2 \mu^5 \lambda, t_1 = \mu^6 - 2\mu\lambda^5, t_2 = 5\lambda^4 \mu^2, t_3 = 5\lambda^2 \mu^4$
of $\OO(-1) \oplus \OO(6)$. $\Box$
\end{prop}

Given a section $s$ of $HM$, we can associate to it the vector $(c_0,
c_1, c_2, c_3)$ representing the coordinates of $s_{L_{00}}$ with
respect to the basis $t_0, t_1, t_2, t_3$.  The coordinates $c_i$ are
then homogeneous coordinates on the moduli space $\PP^3$ of abelian
surfaces in $\PP^4$.  We will call these coordinates \textit{BHM coordinates.}

The polynomial $c_0 t_0 + c_1 t_1 + c_2 t_2 + c_3 t_3$ determines the
singularities of $Z(s)$:

\begin{thm}
Let $s$ be a section of $HM$, and let $f = c_0 t_0 + c_1 t_1 + c_2 t_2 +
c_3 t_3$ be its restriction to $L_{00}$.  The degeneracies of $Z(s)$
are determined by the multiplicities of the roots $(\lambda:\mu)$ of $f$:
\end{thm}

\begin{tabular}{cc}
Multiplicities of roots & Degeneracy of $X(s)$ \\
$(1,1,1,1,1,1)$ & smooth \\
$(2,1,1,1,1,1)$ & translation scroll of elliptic normal curve \\
$(3,1,1,1)$ & tangent scroll of elliptic normal curve \\
$(2,2,1,1)$ & union of five quadrics \\
$(2,2,2)$ & doubled elliptic quintic scroll \\
$(4,2)$ & union of five double planes
\end{tabular}
\newline

The automorphisms $\mu$, $\nu$ and $\delta$ of $\PP^4$ induce
automorphisms, also denoted $\mu$, $\nu$ and $\delta$, of the moduli
space $\PP^3$.  In BHM coordinates, they have the following
form:

\begin{equation}
\mu = \begin{pmatrix}   & 1 &  & \\ 1 &  &  & \\  &  &  & 1 \\  &
& 1 &  \end{pmatrix}, \nu = \begin{pmatrix} \epsilon^4 & & & \\ &
\epsilon & & \\ & & \epsilon^3 & \\ & & & \epsilon^2 \end{pmatrix},
\delta = \frac{1}{\sqrt{5}} \begin{pmatrix} -1 & 1 & \eta' & -\eta \\ 1
& -1 & -\eta & \eta' \\ \eta' & -\eta & 1 & -1 \\ -\eta & \eta' & -1 &
1 \end{pmatrix}
\end{equation}

\noindent where $\eta = \epsilon + \epsilon^4$ and $\eta' = \epsilon^2
+ \epsilon^3$.

\subsection{A pencil of abelian surfaces}

One can ask if the abelian surfaces corresponding to fixed points of
these automorphisms have any interesting properties; this is how we
found the threefold $X$.  Let us find the fixed points of the
automorphism $\mu$; these correspond to the eigenspaces of the matrix

\[
\begin{pmatrix}  & 1 &  &  \\ 1 &  &  &  \\  &  &   & 1 \\  &  &
1 &  \end{pmatrix}.
\]

We find the eigenspace $V_{-1}$ spanned by the vectors $(1,-1,0,0)$
and $(0,0,1,-1)$ and the eigenspace $V_1$ spanned by
$(1,1,0,0)$ and $(0,0,1,1)$  These correspond to lines in the
moduli space $\PP^3$, also denoted $V_{-1}$ and $V_1$.  Our threefold
$X$ will be derived from the threefolds swept out by the abelian surfaces in
$V_{-1}$ and $V_1$.

As in \cite{bib:BM}, we will use the Clebsch diagonal cubic

\begin{equation}
X_3 = \{ c \in \PP^3 : c_0^2 c_3 + c_1^2 c_2 - c_0 c_2^2 - c_1 c_3^2 = 0 \}.
\end{equation}

$X_3$ is the image of the $SL(2,\FF_5)$-equivariant rational map $p: \PP^2
\rightarrow \PP^3$ sending $(y_1 : y_2 : y_3)$ to

\begin{equation}
x = (y_1 y_3^2 - y_2^3 : y_3^3 - y_1 y_2^2 : y_2^2 y_3 - y_2 y_1^2 :
y_3 y_1^2 - y_2 y_3^2).
\end{equation}

The rational map $p$ is undefined at the points $(1:0:0)$ and
$(1:\epsilon^k:\epsilon^{-k})$; these six points correspond to the six
exceptional divisors when we exhibit $X_3$ as $\PP^2$ blown up in six
points.  

Recall the configuration of 27 lines on the cubic surface:  the cubic
surface is isomorphic to $\PP^2$ blown up in six points $p_1, p_2,
\dots, p_6$.  The lines
$E_m$, $m = 1, 2, \dots, 6$ are the exceptional divisors.  The lines
$F_{mn}$, $1 \leq m < n \leq 6$ are the proper transforms of the lines
through $p_m$ and $p_n$.  The lines $G_n$, $n = 1, 2, \dots, 6$ are
the proper transforms of the conics through the five points other than $p_n$.

To locate the six exceptional divisors in $X_3$, temporarily
dehomogenize $y$ by setting $y_1 = 1$.  Also set $y_2 = 0$ and
consider what happens when we let $y_3$ approach 0; we see that
$p((1:0:0))$ approaches the point $(0:0:1:0)$.  Now set $y_3 = 0$ and
consider what happens when $y_2$ approaches 0; $p((1:0:0))$ approaches
the point $(0:0:0:1)$.  Hence the line $E_0$ in $X_3$ is spanned
by the points $(0:0:0:1)$ and $(0:0:1:0)$.  Repeating the same local
analysis at the other five points, we find that the lines $E_k$, $k =
1, 2, \dots, 5$ are spanned by the points
$(-3\epsilon^{2k}:-2\epsilon^k:1:-\epsilon^{-2k})$ and
$(2\epsilon^{-k}:3\epsilon^{-2k}:\epsilon^{2k}:-1)$.

\begin{lem}
The line $V_{-1}$ is exactly the line $F_{56}$.
\end{lem}

\begin{proof} One checks that the line $V_{-1}$ intersects
the $E$-lines $E_5$, $E_6$ and no others.
\end{proof}

Let $F$ be the threefold swept out by the abelian surfaces
parametrized by $F_{56} = V_{-1}$.  Let $G$ be the threefold swept out
by $V_1$.  We can determine what all the fibres of $F$ are:  

\begin{prop}
The singular fibres of $F$ are as follows:  there are two fibres of
type (2,2,1,1) corresponding to unions of five quadrics and two fibres
of type (2,2,2) corresponding to doubled elliptic quintic scrolls.

The singular fibres of $G$ are as follows:  there are two fibres of
type (2,2,1,1) corresponding to unions of five quadrics and two fibres
of type (3,1,1,1) corresponding to tangent scrolls of elliptic normal curves.
\end{prop}

\begin{proof} For $(a:b) \in \PP^1$ and $(\lambda:\mu) \in
\PP^1$, consider the correspondence in $\PP^1 \times \PP^1$ defined by
the condition that $(\lambda:\mu)$ be a root of the polynomial $a(t_0
- t_1) + b(t_2 - t_3)$.  This is a $(1,6)$ correspondence and thus has
genus zero.  Projection to the first copy of $\PP^1$ is a map between
rational curves of degree 6.  By the Riemann-Hurwitz theorem, there
are 10 ramification points (counted with multiplicity).

The line $F_{56}$
intersects $E_5$ and $E_6$ at $(2:-2:4:4)$ and $(0:0:1:-1)$
respectively.  The corresponding sextic polynomials in $\lambda, \mu$
have multiplicities $(2,2,1,1)$, so the corresponding fibres are
unions of five quadrics each.

Four branch points have already been accounted for.  We found that the
line $F_{56}$ intersects the curve $C_6$ in the points $(5:-5:2i+1:-2i-1)$ and
$(5:-5:-2i+1:2i-1)$.  These points correspond to polynomials of type
$(2,2,2)$.  They account for the other six ramification points.  We
have therefore found all the singular fibres of $F$.

For $G$, we look at the correspondence on $\PP^1 \times \PP^1$ where
$(\lambda:\mu)$ is a root of $a(t_0 + t_1) + b(t_2 + t_3)$.  This
time, the correspondence has two horizontal components because $(1:i)$
and $(1:-i)$ are common roots of $t_0 + t_1$ and $t_2 + t_3$.  The
remaining part of the correspondence is a $(1,4)$ curve that maps
surjectively onto the first copy of $\PP^1$, so there are 6
ramification points.

By inspection, the polynomial $a(t_0 + t_1) + b(t_2 + t_3)$ has
repeated roots when $(a:b)$ = $(1:0)$, $(0:1)$, $(1:\frac{3+4i}{5})$
and $(1:\frac{3-4i}{5})$ of type $(2,2,1,1)$, $(2,2,1,1)$, $(3,1,1,1)$
and $(3,1,1,1)$ respectively.  Eliminating the common roots $(1:\pm
i)$, the remaining roots have multiplicity $(2,2)$, $(2,2)$,
$(2,1,1)$ and $(2,1,1)$ respectively.  This accounts for the 6
ramification points.
\end{proof}

Once we have the defining equation for $F$, simple calculations will show that the quadric surface
$T_0 = \{ x: x_0 = x_1 x_4
- x_2 x_3 = 0 \}$ and its translates $T_i = \sigma^i T_0$ are
contained in $F$; the $T_i$ are $\tau$-invariant.  Put $U_0 = \delta T_0$; 
we then have $U_0 = \{ x:
\Sigma_i x_i = \Sigma_{i \neq j} x_i x_j = 0 \}$.  The quadric surface
$U_0$ and its translates $U_i = \tau^i W_0$ are contained in $F$ as
well; the $U_i$ are $\sigma$-invariant.  Since the surfaces $\cup_i
T_i$ and $\cup_i U_i$ are of degree 10 and invariant under $H_5$, they
must be two of the singular fibres in the abelian surface fibration of
$F$.

Similarly, one checks that the quadric surface $Q_0 = \{ z: z_0 = z_1
z_4 + z_2 z_3 = 0 \}$ and its translates $Q_i = \sigma^i Q_0$ are in
$G$, and that the quadric surface $R_0 = \delta Q_0 = \{ z: \Sigma_i
z_i = \Sigma_{i \neq j} z_i z_j = 0 \}$ and its translates $Q_i =
\tau^i Q_0$ are in $G$.  These two unions of five quadrics are two of
the singular fibres in the fibration of $G$.

We can identify the other fibres of $F$; first we need to identify certain 
elliptic curves in $\PP^4$.  Aure, Decker, Hulek, Popescu and Ranestad 
\cite{bib:ADHPR} have shown that the set
of $H_5$-equivariant elliptic normal curves in $\PP^4$ is parametrized
by $\PP^1$:  to the point $(\lambda:\mu)$ corresponds the elliptic
normal curve $E_{(\lambda:\mu)}$ defined by the set of equations

\begin{equation}
q_i^{(\lambda:\mu)}(x) = -\lambda\mu x_i^2 - \mu^2 x_{i+1} x_{i+4} +
\lambda^2 x_{i+2} x_{i+3}.
\end{equation}

By inspection, we see that the curves corresponding to $(1:i)$ and
$(1:-i)$ are in $F$ and $G$.  Denote these curves by $E_1$ and $E_2$
respectively.  Proposition 4.3 in \cite{bib:ADHPR} shows that
the elliptic quintic scrolls $Q_1$ and $Q_2$ contain $E_1$ and $E_2$ respectively, where $Q_1$
is defined by the equations

\begin{equation}
x_i^3 + x_i x_{i+1} x_{i+4} + x_i x_{i+2} x_{i+3} - i (x_{i+1}^2
x_{i+3} + x_{i+2} x_{i+4} + x_{i+1} x_{i+2}^2 + x_{i+3}^2 x_{i+4}) = 0
\end{equation}

\noindent and $Q_2$ is defined by replacing $i$ above with its complex conjugate
$-i$.  Simple calculations show that $Q_1$ and $Q_2$ are contained in
$F$.  Being $H_5$-invariant, they must be the elliptic quintic scroll
fibers of $F$.


\section{The threefold $\XT$}

\subsection{Definition of $\XT$}

Although we have defined the Horrocks-Mumford bundle only over $\CC$, the
pencils of abelian surfaces it defines are quintics in $\PP^4$, and
the quintics we are interested in have integer coefficients and can
thus be studied over arbitrary fields $k$.

In particular, we have the Horrocks-Mumford quintics $F$ and $G$ defined by the lines
$V_{-1}$ and $V_1$ in $\PP^3$.  In \cite{bib:Manolache}, Manolache
was able to determine the equation of the Horrocks-Mumford quintic
in terms of the BHM-coordinates of the parametrizing line:

\begin{thm}
Suppose the Horrocks-Mumford quintic $X$ is determined by the line
passing through the points $(a_0:a_1:a_2:a_3)$ and
$(b_0:b_1:b_2:b_3)$.  Then the defining equation for $X$ is
\begin{equation}
\begin{aligned}
& \quad 25 (a_3 b_4 - a_4 b_3) (\Sigma x_0 x_1 x_2 x_3 x_4) \\
&+ 5 (a_2 b_3 - a_3 b_2) (\Sigma x_0 x_2^2 x_3^2) \\
&- 5 (a_1 b_3 - a_3 b_1) (\Sigma x_0^3 x_2 x_3) \\
&+ 5 (a_2 b_4 - a_4 b_2) (\Sigma x_0^3 x_1 x_4) \\
&- 5 (a_1 b_4 - a_4 b_1) (\Sigma x_0 x_1^2 x_4^2) \\
&+ (a_1 b_2 - a_2 b_1) (\Sigma (x_0^5 - x_0 x_1 x_2 x_3 x_4) \\
\end{aligned}
\end{equation}
where the sums are taken over cyclic permutations of the indices.
\end{thm}

An easy calculation shows that $F$ is defined by the equation

\begin{equation}
\Sigma_i (x_i^3 x_{i+1} x_{i+4} + x_i^3 x_{i+2} x_{i+3} - x_i x_{i+1}^2
x_{i+4}^2 - x_i x_{i+2}^2 x_{i+3}^2) = 0
\end{equation}

\noindent and that $G$ is defined by the equation

\begin{equation}
\Sigma_i ( z_i^3 z_{i+1} z_{i+4} - z_i^3 z_{i+2} z_{i+3} - z_i z_{i+1}^2
z_{i+4}^2 + z_i z_{i+2}^2 z_{i+3}^2) = 0.
\end{equation}

\noindent where the summations are taken over cyclic permutations of the
indices.  As before, $F$ and $G$ are both invariant under the action of the
matrices $\sigma$ and $\tau$.

Consider now the complete intersection threefold $X$ in $\PP^4(x) \times
\PP^4(z)$ given by the matrix equation

\begin{equation}
M(x) z = 0,
\end{equation}

\noindent where

\begin{equation}
M(x) = \begin{pmatrix}  & -x_3 & x_1 & x_4 & -x_2 \\ -x_3 &  & -x_4 &
  x_2 & x_0 \\ x_1 & -x_4 &  & -x_0 & x_3 \\ x_4 & x_2 & -x_0 & 
  & -x_1 \\ -x_2 & x_0 & x_3 & -x_1 &  \end{pmatrix}.
\end{equation} 

Note that this is equivalent to the matrix equation

\begin{equation}
L(z) x = 0,
\end{equation}

\noindent where

\begin{equation}
L(z) = \begin{pmatrix}  & z_2 & -z_4 & -z_1 & z_3 \\ z_4 &  & z_3
  & -z_0 & -z_2 \\ -z_3 & z_0 &  & z_4 & -z_1 \\ -z_2 & -z_4 & z_1 &
   & z_0 \\ z_1 & -z_3 & -z_0 & z_2 &  \end{pmatrix}.
\end{equation}

Note also that $\det M(x)$ and $\det L(z)$ give us the equations for
$F$ and $G$ respectively (up to a factor of 2).  Hence the projections $\pi_1$ and $\pi_2$ of
$X$ onto each factor give us $F$ and $G$.  In \cite{bib:Moore},
Moore first considered the matrices

\begin{equation}
M(x,y) = \begin{pmatrix} x_0 y_0 & x_3 y_2 & x_1 y_4 & x_4 y_1 & x_2
  y_3 \\
x_3 y_3 & x_1 y_0 & x_4 y_2 & x_2 y_4 & x_0 y_1 \\
x_1 y_1 & x_4 y_3 & x_2 y_0 & x_0 y_2 & x_3 y_4 \\
x_4 y_4 & x_2 y_1 & x_0 y_3 & x_3 y_0 & x_1 y_2 \\
x_2 y_2 & x_0 y_4 & x_3 y_1 & x_1 y_3 & x_4 y_0 \end{pmatrix},
\end{equation}

\begin{equation}
L(z,y) = \begin{pmatrix} z_0 y_0 & z_2 y_4 & z_4 y_3 & z_1 y_2 & z_3
  y_1 \\
z_4 y_1 & z_1 y_0 & z_3 y_4 & z_0 y_3 & z_2 y_2 \\
z_3 y_2 & z_0 y_1 & z_2 y_0 & z_4 y_4 & z_1 y_3 \\
z_2 y_3 & z_4 y_2 & z_1 y_1 & z_3 y_0 & z_0 y_4 \\
z_1 y_4 & z_3 y_3 & z_0 y_2 & z_2 y_1 & z_4 y_0 \end{pmatrix}.
\end{equation}

\noindent For a generic choice of $y \in \PP^4$, the threefold
determined by $\det M(x,y) = 0$ and the threefold determined by $\det
L(z,y) = 0$ are both Horrocks-Mumford quintics with the expected 100
nodes.  For our threefolds $F$ and $G$, we have taken $y = (0:1:-1:-1:1)$.

\subsection{Singularities of $X$}

We need to know the singularities of $X$:

\begin{prop}
Over a field of characteristic not equal to 2 or 5, $X$ has 60
singular points.  The 60 singular
points of $X$ are all ordinary double points (nodes).\label{theorem:main}
\end{prop}

\begin{rem} In characteristic 0, Gross and Popescu \cite{bib:GP} have studied
the two-parameter family of threefolds $X_y$ given by perfoming
the construction above for $y$ a generic point of the plane

\begin{equation*}
\PP^2_{+} = \{ y : y_1 - y_4 = y_2 - y_3 = 0 \}.
\end{equation*}

Thus $X_y$ is a common partial resolution of the quintic
threefolds $F_y$ and $G_y$.  Our threefold $X$ is thus a special
member of this family.  Over $\CC$, Gross and Popescu proved
the following statements for a generic choice of $y$:

\begin{enumerate}

\item $F_y$ is singular along the union of two elliptic curves
  $D_{1,y}$ and $D_{2,y}$.  These curves are the base curves of the
  elliptic quintic scrolls $Q_{1,y}$ and $Q_{2,y}$ appearing in the
  abelian surface fibration of $F_y$, and they intersect
  in the 25 points comprising the $H_5$ orbit of $y$.  These 25 points
  comprise the base locus of the pencil of abelian surfaces.

\item $G_y$ is singular along a disjoint union of two elliptic curves
  $E_{1,y}$ and $E_{2,y}$.

\item $X_y$ is singular along 50 nodes, 2 of which lie over each
  point of the Heisenberg orbit of $y$ when projecting via the map
  $\pi_1$.  We call these 50 nodes {\em regular nodes.}

\end{enumerate}
\end{rem}

\noindent Essentially our threefold $X$ is a special member of the family
$X_y$ where $y = (0:1:-1:-1:1)$.  Over $\CC$, it is easy to check
that $X$ has the 50 expected nodes and 10 others.  Over fields of
arbitrary characteristic, we resort to brute-force computation.  A
basic trick we use is to compute Gr\"obner bases of ideals over $\ZZ$ in
order to obtain results valid in fields of unknown characteristic.

We will begin the proof of the proposition after we dispense with some preliminary
results.  If $x$ is a point in $\PP^4(x)$, we say that $x$ is a rank
$n$ point if rank $M(x)$ = $n$.  Similarly, if $z$ is a point in
$\PP^4(z)$, we say that $z$ is a rank $n$ point if rank $L(z)$ = $n$.

\begin{lem}
If $x$ is a point of $F$, then $x$ is a rank 4 point or a rank 3
point.
\end{lem}

\begin{proof}  If $x$ is a point of $F$, then $\det M(x) =
0$.  Hence the rank of $\det M(x)$ is at most 4.

The 3 by 3 matrices

\[
\begin{pmatrix}
 & -x_3 & x_1 \\ -x_3 &  & -x_4 \\ x_1 & -x_4 & 
\end{pmatrix}
\begin{pmatrix}
 & -x_4 & x_2 \\ -x_4 &  & -x_0 \\ x_2 & -x_0 & 
\end{pmatrix}
\begin{pmatrix}
 & -x_0 & x_3 \\ -x_0 &  & -x_1 \\ x_3 & -x_1 & 
\end{pmatrix}
\]
\[
\begin{pmatrix}
 & x_1 & -x_2 \\ x_1 &  & x_3 \\ -x_2 & x_3 &
\end{pmatrix}
\begin{pmatrix}
 & -x_3 & -x_2 \\ -x_3 &  & x_0 \\ -x_2 & x_0 & 
\end{pmatrix}
\begin{pmatrix}
 & x_1 & x_4 \\ x_1 &  & -x_0 \\ x_4 & -x_0 & 
\end{pmatrix}
\]
\[
\begin{pmatrix}
 & x_2 & x_0 \\ x_2 &  & -x_1 \\ x_0 & -x_1 &
\end{pmatrix}
\begin{pmatrix}
 & x_1 & -x_2 \\ x_1 &  & x_3 \\ -x_2 & x_3 &
\end{pmatrix}
\begin{pmatrix}
 & -x_3 & x_4 \\ -x_3 &  & x_2 \\ x_4 & x_2 & 
\end{pmatrix}
\]
\[
\begin{pmatrix}
 & -x_4 & x_0 \\ -x_4 &  & x_3 \\ x_0 & x_3 & 
\end{pmatrix}
\]

\noindent are all 3 by 3 minors of $M(x)$.  Their determinants give us $\pm 2
x_i x_j x_k$ for all subsets $\{ i,j,k \}$ of $\{ 0, 1, 2, 3, 4 \}$.
Suppose these determinants are all zero.  Then three of the $x_i$ must
be zero.  

Without loss of generality, we can assume either that $x_0, x_1, x_2$
are 0 or that $x_0, x_2, x_4$ are zero.  Suppose the former holds.  One sees
then that the matrix

\[\begin{pmatrix}
 & -x_3 &  & x_4 &  \\ -x_3 &  & -x_4 &  &  \\  & -x_4 &  & 
& x_3 \\ x_4 &  &  &  &  \\  &  & x_3 &  & 
\end{pmatrix}
\]

\noindent has at least three linearly independent columns.  A similar result
holds of $x_0, x_2, x_4$ are zero.  Hence the rank of
$M(x)$ is at least 3 for all $x$.
\end{proof}

\begin{lem}
If $z$ is a point of $G$, then $z$ is a rank 4 point or a rank 3
point.
\end{lem}

\begin{proof}  Again, if $z$ is a point of $G$, then $\det
L(z) = 0$.  Hence the rank of $\det L(z)$ is at most 4.

The polynomials $\pm (z_i z_{i+1}^2 + z_{i+2}^2 z_{i+3})$ and $\pm (z_i^2 z_{i+3} +
z_{i+2}^2 z_{i+4})$ are all determinants of 3 by 3 minors of $L(z)$.
Suppose that these minors are all zero.  Without loss of generality,
assume $z_0 = 1$.  

Since $z_0 z_1^2 + z_2^2 z_3 = 0$, we have $z_1^2 = -z_2^2 z_3$.  From
$z_4 z_0^2 + z_1^2 z_2 = 0$, we get $z_4 = -z_1^2 z_2 = z_2^3 z_3$.
From $z_2 z_3^2 + z_4^2 z_0 = 0$, we get $z_2 z_3^2 (1 + z_2^5) = 0$.

On the other hand, we also have $z_4^2 z_2 + z_1^2 z_3 = 0$.  From
this we see that $z_2^7 z_3^2 - z_2^2 z_3^2 = 0$, or $(z_2^5 - 1)
z_2^2 z_3^2 = 0$.  Thus we must have either $z_2 = 0 $ or $z_3 = 0$.
In either case, we have $z = (1:0:0:0:0)$, which is a rank 4 point.
Hence there are no points of rank lower than 3.  
\end{proof}

\begin{lem}  Points of $F$ or $G$ of rank less than 4 are singular
  points.
\end{lem}

\begin{proof}  More generally, let $A$ be an $n$ by $n$ matrix
whose $ij$ entry is the indeterminate $x_{ij}$.  Using the Laplace
expansion formula, we see that $\frac{\partial \det(A)}{\partial x_{ij}}$
is equal to $(-1)^{i+j} \det A_{ij}$, where $A_{ij}$ is the $ij$ minor
of $A$.

Now suppose that the $x_{ij}$ are functions of some other
indeterminates $y_k$.  By the chain rule, $\frac{\partial
  \det(A)}{\partial y_k}$ = $\sum^n_{i = 1} \sum^n_{j = 1} (-1)^{i+j}
\det A_{ij} \frac{\partial x_{ij}}{\partial y_k}$.  If $y$ is a point
of rank less than $n$, then the determinants of all the $(n-1)$ by
$(n-1)$ minors are zero, and hence $\frac{\partial \det(A)}{\partial
  y_k}$ is zero for all $k$.
\end{proof}

\begin{lem}  The projections $\pi_1: X \rightarrow F$ and $\pi_2:
  X \rightarrow G$ are isomorphisms outside the rank 3 loci of
  $F$ and $G$ respectively.
\end{lem}

\begin{proof}  If $x$ is a rank 4 point of $\PP^4(x)$, this
means that the kernel of $M(x)$ is 1-dimensional.  Hence the kernel of
$M(x)$ defines a unique point $z$ in $\PP^4(z)$.  We thus have a
regular map $f_1 : F - F_4 \rightarrow \pi^{-1}(F - F_4)$.

Now given a point $(x,z)$ in $X$ such that $x$ is a rank 4 point,
the coordinates of $z$ are given by the determinants of 4 by 4 minors
of $M(x)$.  Hence the map $\pi_1: \pi^{-1}(F - F_4) \rightarrow (F -
F_4)$ is also regular, proving the result for $F$.  A similar result
holds for $G$.
\end{proof}

\begin{lem}  Suppose the characteristic of the base field is not 2.
  Then the rank 3 locus of $G$ is the union $E$ of $E_1$ and $E_2$, where $E_1$ is cut out by the
  polynomials $i z_i^2 + z_{i+1} z_{i+4} + z_{i+2} z_{i+3}$ and $E_2$
  is cut out by the polynomials $-i z_i^2 + z_{i+1} z_{i+4} + z_{i+2}
  z_{i+3}$.
\end{lem}

\begin{proof} This is proven by using Macaulay2.  One shows
that the radicals of the ideals defining the two varieties have
identical Gr\"obner bases, up to constant factors of 2 and 5.
\end{proof}

\begin{question} Is there a geometric explanation for this
result?
\end{question}

\begin{lem}  If the characteristic of the base field is not 2 or 5, $E_1$ and $E_2$ are elliptic normal curves.
\end{lem}

\begin{proof} This was proven by Fisher in Chapter 4 of
\cite{bib:Fisher}.  Assuming the characteristic of the base field is
not 5, Fisher constructs the universal curve $\mathcal{X}(5) \subset
X(n) \times \PP^4$ as the closure of the scheme defined by the 4 by 4
Pfaffians of the 5 by 5 matrix 

\[
\begin{pmatrix}
         & -a_1 x_1 & -a_2 x_2 &  a_2 x_3  &  a_1 x_4 \\
 a_1 x_1 &          & -a_1 x_3 & -a_2 x_4  &  a_2 x_0 \\
 a_2 x_2 &  a_1 x_3 &          & -a_1 x_0  & -a_2 x_1 \\
-a_2 x_3 &  a_2 x_4 &  a_1 x_0 &           & -a_1 x_2 \\
-a_1 x_4 & -a_2 x_0 &  a_2 x_1 &  a_1 x_2  &          \\
\end{pmatrix}
\]

\noindent where $a = (0:a_1:a_2:-a_2:-a_1)$.  One checks that $E_1$ and $E_2$ are the
curves obtained when $a = (0:1:-i:i:-1)$ and $(0:1:i:-i:1)$ respectively.  By considering the $SL_2(\ZZ/ 5\ZZ)$ action on
$\mathcal{X}(5)$, Fisher shows that the fibers are smooth elliptic
normal curves when $a_1 a_2 (a_1^{10} - 11 a_1^5 a_2^5 - a_2^{10})
\neq 0$, in analogy with the case where the base field is $\CC$.
\end{proof}

Note that $S_1 = \pi_2^{-1}(E_1)$ and $S_2 = \pi_2^{-1}(E_2)$ are elliptic ruled
surfaces in $X$.

We can now classify the singularities of $X.$

\begin{proof}[Proof of proposition \ref{theorem:main}.]
A point $(x,z)$ of $X$ is
a singular point if and only if the kernel of the matrix
$\begin{pmatrix} L(z) & M(x) \end{pmatrix}$ has dimension at least 6,
i.e. if the rank is at most 4.  This is equivalent to saying that the
rank of the transposed matrix $\begin{pmatrix} L^T(z) \\ M(x)
\end{pmatrix}$ is at most 4, which is equivalent to saying that the
kernel of $\begin{pmatrix} L^T(z) \\ M(x) \end{pmatrix}$ has dimension
at least 1, i.e. the kernel of $L^T(z)$ and the kernel of $M(x)$ have
nontrivial intersection.

\begin{case1}  $x$ is a rank 4 point.  
\end{case1}
Since $(x,z)$ is a
point of $X$, $z$ is in the kernel of $M(x)$.  Therefore $z$ must
span the space $\ker L^T(z) \cap M(x)$.  Hence $L^T(z) z = 0$.  Simple
algebra shows that the only singular points $(x,z)$ in this case are
$((1:0:0:0:0),(1:0:0:0:0))$ and $((1:1:1:1:1),(1:1:1:1:1))$ and their
orbits under $(\sigma, \sigma)$ and $(\tau^3, \tau)$.  We will call
these nodes in $X$ $\sigma-nodes$ and $\tau-nodes$ respectively.
We will also use the same nomenclature for the images of these nodes
in $F$ and $G$.

\begin{case2}
$x$ is a rank 3 point.
\end{case2}
  We will suppose that $z$
must then be a rank 3 point.  Using Macaulay2, we show that if $(x,z)$
is a singular point of $X$ and $x$ is a rank 3 point of $F$,
then some $x_i$ or some $z_i$ is zero.  (The code can be found in the Appendix.)

Suppose that some $x_i$ is zero.  Without loss of
generality, suppose $x_0$ is zero.  The polynomials $(x_1 x_2 +
x_4^2)^2$, $(x_2 x_4 + x_3^2)^2$, $(x_1 x_3 + x_2^2)^2$, $(x_1 x_4 +
x_3^2)^2$ and $(x_2 x_3 - x_1 x_4)^2$ are all $4 \times 4$ minors of
$M(x)$ when $x_0$ is zero.

If $x_1$ is also zero, we quickly see that all $x_i$ are zero, a
contradiction.  So suppose $x_1 = 1$.  We quickly see that $x_4 =
\epsilon^k$, $x_2 = -\epsilon^{2k}$ and $x_3 = -\epsilon^{4k}$.
Solving for $z$ gives us the singular points

\[
(x,z) = ((0:1:-\epsilon^k:-\epsilon^{2k}:\epsilon^{3k}), (0:1:\pm
i \epsilon^{2k}: \mp i \epsilon^{4k}: -\epsilon^{k})),
\]

\noindent whose orbits under $(\sigma, \sigma)$ will be the 50 regular nodes of $X$.

Now suppose instead that some $z_i$ is zero.  Without loss of generality, we may assume
$z_0 = 0.$.  If we assume that $z$ is a rank 3 point of $G$, then $z$
is a point of $E = E_1 \cup E_2$ with $z_0 = 0$.  We have equations
defining $E$, and using these equations we recover the 50 points above.

If $z$ is a rank 4 point, then $z$ must itself be a rank 4 node of
$G$, since the map $\pi_2: X \rightarrow G$ is an isomorphism on
the rank 4 locus.  Again using Macaulay2, we can find generators of an
ideal that defines the set of singular points on $G$ with $z_0 = 0$.
We eventually find one of two things:  either $z$ is a $\sigma$-node, all of
whose corresponding $x$ are not rank 3 points, or $z$ is of the form
$(0: 1: \pm i \epsilon^{k}: \mp i \epsilon^{2k}: -\epsilon^{3k})$,
which are not rank 4 points of $G$.  So we have found no new nodes.

We have assumed throughout that our nodes have $x_0 = 0$ or $z_0 =
0$; taking the orbits of $(x,z)$ under $(\sigma, \sigma)$, we get all possible values of
$z$ for which there is a node $(x,z)$ with $x$ and $z$ rank 3 points.

There are 60 singular points on $X$:  the $\sigma$-nodes; the
$\tau$-nodes; and the 50 regular nodes given by
$((0:1:-1:-1:1),(0:1:\pm i:-\pm i:-1))$ and their orbits under $(\sigma,\sigma)$ and $(\tau^3,\tau)$.

To prove that the 60 singular points are indeed nodes, we need to construct
local coordinates around each singular point.

$((1:0:0:0:0),(1:0:0:0:0))$:  Dehomogenize by fixing $x_0 = 1$ and
$z_0 = 1$.  $X$ is then defined by the five equations

\begin{equation}
\begin{aligned}
-x_3 z_1 + x_1 z_2 + x_4 z_3 - x_2 z_4 &= 0, \\
-x_3 - x_4 z_2 + x_2 z_3 + z_4 &= 0, \\
x_1 - x_4 z_1 - z_3 + x_3 z_4 &= 0, \\
x_4 + x_2 z_1 - z_2 - x_1 z_4 &= 0, \\
-x_2 + z_1 + x_3 z_2 - x_1 z_3 &= 0. \\
\end{aligned}
\end{equation}

Let $I$ be the ideal generated by these five polynomials.
  
Notice that $z_4 = x_3$ plus terms of higher
order after we complete the coordinate ring at the ideal
$(x_1,\dots,x_4,z_1,\dots,z_4)$.  Similarly $z_3 = x_1 + h.o.t.$, $z_2 = x_4 + h.o.t.$ and $z_1
= x_2 + h.o.t.$.  

Replacing the $z$'s in this manner, we see that $X$ is defined by the single equation

\begin{equation}
-x_3 x_2 + x_1 x_4 + x_4 x_1 - x_2 x_3 + h.o.t. = 0 \label{equation:quadric}
\end{equation}

\noindent in the variables $x_1, x_2, x_3, x_4$.  Let $J$ be the ideal
generated by this formal power series.  We see that the completed local ring
$\kbar[[x_1,\dots,x_4,z_1,\dots,z_4]]/I$ is isomorphic to the completed local ring $\kbar[[x_1,\dots,x_4]]/J$.

The second-order piece of equation (\ref{equation:quadric}) reads

$x^{T} B x = 0$ where $B = \begin{pmatrix}  &  &  & 1 \\  &  & -1
  &  \\  & -1 &  &  \\ 1 &  &  &  \end{pmatrix}$.  The
  determinant of $B$ is 1.  Therefore, we can find a matrix $C$ such
  that $B = C^T C$.  Put $y = Cx$; we now have an isomorphism

\[
\kbar[[x_1, \dots, x_4]]/J \longrightarrow \kbar[[y_1, \dots, y_4]]/L
\]

\noindent where $L$ is generated by some formal power series of the form
$y_1^2 + y_2^2 + y_3^2 + y_4^2 + h.o.t.$ and $Cx \rightarrow y$.
 
Finally, we show that there is an isomorphism

\[
\kbar[[y_1 \dots y_4]]/L \rightarrow \kbar[[t_1, \dots, t_4]]/M
\]

\noindent where $M$ is generated by $t_1^2 + t_2^2 + t_3^2 + t_4^2$.  Put $t_i =
y_i + P_{i2}(y) + P_{i3}(y) + \dots$, where $P_{ik}$ is a monomial of
order $k$ in the $y$'s.  $L$ is generated by $y_1^2 + y_2^2 + y_3^2 +
y_4^2 + L_3(y) + \dots$.  We can find $P_{i2}$ such that $\sum_i 2 y_i
P_{i2} = L_3(y)$, and we can solve for higher $P_{ik}$ recursively.
The resulting expressions for the $t_i$ define the isomorphism.  At
last we see that our singular point is a node.

$((1:1:1:1:1),(1:1:1:1:1))$:  Dehomogenize by setting $x_0 = 1$ and
  $z_0 = 1$.  Now put $x_i = u_i + 1$ and $y_i = w_i + 1$ for $i =
  1,2,3,4$.  The defining equations for $X$ are now

\begin{equation}
\begin{aligned}
u_1 - u_2 - u_3 + u_4 - w_1 + w_2 + w_3 - w_4 + (-u_3 w_1 + u_1 w_2 +
u_4 w_3 - u_2 w_4) &= 0, \\
u_2 - u_3 - u_4 - w_2 + w_3 + w_4 + (-u_4 w_2 + u_2 w_3) &= 0, \\
u_1 + u_3 - u_4 - w_1 - w_3 + w_4 + (-u_4 w_1 + u_3 w_4) &= 0, \\
-u_1 + u_2 + u_4 + w_1 - w_2 - w_4 + (u_2 w_1 - u_1 w_4) &= 0, \\
-u_1 - u_2 + u_3 + w_1 + w_2 - w_3 + (u_3 w_2 - u_1 w_3) &= 0. \\
\end{aligned}
\end{equation}

Adding them up, we get 

\begin{equation}
\begin{split}
&u_1 w_2 - u_1 w_3 - u_1 w_4 + u_2 w_1 + u_2 w_3 - u_2 w_4 \\
&- u_3 w_1 + u_3 w_2 + u_3 w_4 - u_4 w_1 - u_4 w_2 + u_4 w_3 = 0. \\
\end{split}
\end{equation}

Again we can eliminate the $w_i$; we then have

\begin{equation}
\begin{split}
&u_1 (u_2 - u_3 - u_4) + u_2 (u_1 + u_3 - u_4) \\
&+ u_3 (-u_1 + u_2 + u_4) + u_4 (-u_1 - u_2 - u_3) + h.o.t. = 0.\label{equation:quadric2}
\end{split}
\end{equation}

\noindent after passing to the completion of the coordinate ring.

The second-order piece of (\ref{equation:quadric2}) can be written as

\[
2 u^T B u = 0
\]

\noindent with $B = \begin{pmatrix}  & 1 & -1 & -1 \\ 1 &  & 1 & -1 \\ -1 & 1
  &  & 1 \\ -1 & -1 & 1 &  \end{pmatrix}$.  The determinant of $B$
  is 5.  Therefore, $((1:1:1:1:1),(1:1:1:1:1))$ is a node.

$((1:-1:-1:1:0),(1:i:-i:-1:0))$:  Again dehomogenize by setting $x_0 =
  z_0 = 1$.  Put $x = u + (-1,-1,1,0), z = w + (i,-i,1,0)$.  In the
  $(u,w)$ coordinates, the defining equations for $X$ are

\begin{equation}
\begin{aligned}
-i u_1 -i u_3 -u_4 - w_1 -w_2 + w_4 + (u_1 w_2 - u_2 w_4 - u_3 w_1 +
 u_4 w_3) &= 0, \\
-u_2 - u_3 + i u_4 - w_3 + w_4 + (u_2 w_3 - u_4 w_2) &= 0, \\
u_1 - i u_4 - w_3 + w_4 + (u_3 w_4 - u_4 w_1) &= 0, \\
i u_2 + u_4 - w_1 - w_2 + w_4 + (-u_1 w_4 + u_2 w_1) &= 0, \\
u_1 - u_2 - i u_3 + w_1 + w_2 + w_3 + (-u_1 w_3 + u_3 w_2) &= 0. \\
\end{aligned}
\end{equation}

From these equations, we can obtain the equation

\begin{equation}
\begin{split}
&-u_1 w_2 + u_2 w_4 + u_3 w_1 - u_4 w_3 + i u_2 w_3 \\
&- i u_4 w_2 - i u_3 w_4 + i u_4 w_1 - u_1 w_4 + u_2 w_1 = 0. \\
\end{split}
\end{equation}

This time we want to eliminate $u_1, u_2, u_3$ and $w_1$.  We have

\begin{equation}
\begin{aligned}
u_1 &= i u_4 + w_3 - w_4 + h.o.t., \\
u_2 &= i u_4 + i w_3 + \frac{2-i}{5} w_4 + h.o.t., \\
u_3 &= (-1-i) w_3 + \frac{3+i}{5} w_4 + h.o.t., \\
w_1 &= -w_2 - w_3 + \frac{6+2i}{5} w_4 + h.o.t. \\
\end{aligned}
\end{equation}

We then get the equation

\begin{equation}
\begin{split}
&-4i u_4 w_2 - (2 + 2i) u_4 w_3 \\
&- \frac{4-12i}{5} u_4 w_4 - (4-2i) w_3 w_4 + \frac{14-2i}{5} w_4^2 +
  h.o.t. = 0. \\
\end{split}
\end{equation}
The second-order part of this equation can be written in matrix form:

\begin{equation}
v^T B v = 0
\end{equation}

\noindent with $v = (u_4, w_2, w_3, w_4)$ and 

\begin{equation*}
B = \begin{pmatrix}  & -2i &
  -1-i & \frac{-2+6i}{5} \\ -2i &  &  &  \\ -1-i &  &  & -2+i \\
  \frac{-2+6i}{5} &  & -2+i & \frac{14-2i}{5} \end{pmatrix}.
\end{equation*}  

\noindent The determinant of this matrix is $-12 + 16i$.  So this
point is a node.  By symmetry, all of our 60 singular points are nodes.
\end{proof}

Note that the surface $S_1$ contains the 25 nodes

\begin{equation*}
(\sigma^i \tau^3j, \sigma^i \tau^j) ((0:1:-1:-1:1),(0:1:-i:i:-1))
\end{equation*}

and that $S_2$ contains the 25 nodes

\begin{equation*}
(\sigma^i \tau^3j, \sigma^i \tau^j) ((0:1:-1:-1:1),(0:1:i:-i:-1).
\end{equation*}

Furthermore, the surfaces $\pi_1^{-1} T_i$ contain the 5
$\sigma$-nodes and the surfaces $\pi_1^{-1} U_i$ contain the 5
$\tau$-nodes.  Hence we can blow up these surfaces in some order to
obtain a projective small resolution $\hat{X}$ of $X$ over $\QQ(\epsilon)$.  However, we cannot
obtain a model for any small resolution over $\QQ$ because the determinant of the
local quadratic equation at the node $((1:1:1:1:1),(1:1:1:1:1))$, which is not a square in $\QQ$.

Assume from now on that the characteristic of our base field is
different from 2 and 5.  Since the singularities of $X$ are 60
nodes, we can blow up the nodes to get a smooth threefold $\XT$.
(Note that the equations defining the union of the 60 points are
defined over $\ZZ$, so the blowup is defined over $\ZZ$.)  Each
exceptional divisor is isomorphic to $\PP^1 \times \PP^1$, and we will
need to know when the rulings on each $\PP^1 \times \PP^1$ are defined
over $\FF_p$:

\begin{prop}
If $p \equiv 1 ( \Mod 20)$, then all 60 nodes are defined over $\FF_p$.
If $p \equiv 11 ( \Mod 20)$, then only the 10 $\sigma$-nodes and
$\tau$-nodes are defined over $\FF_p$.  If
$p \equiv 9,13,17 ( \Mod 20)$, then the $\sigma$-nodes
$((1:1:1:1:1),(1:1:1:1:1))$ and 10 of the regular nodes are defined over $\FF_p$.
Otherwise only the $\sigma$-nodes and $((1:1:1:1:1),(1:1:1:1:1))$ are defined over $\FF_p$.
\end{prop}

\begin{proof}
  Whether or not our nodes are defined over
$\FF_p$ depends on the values of $p$ modulo 4 and 5.  If $p \equiv 1
(\Mod 4)$ and $p \equiv 1 (\Mod 5)$, then both $i$ and $\epsilon$ are in
$\FF_p$, and hence all 60 nodes are defined over $\FF_p$.  If $p \equiv 3 (\Mod 4)$ and $p \equiv 1 (\Mod 5)$, then the $\sigma$-nodes
and $\tau$-nodes are all defined over $\FF_p$.  If $p \equiv 1 (\Mod 4)$ and $p \equiv
2,3,4 (\Mod 5)$, then the $\sigma$-nodes, $((1:1:1:1:1),(1:1:1:1:1))$
and the nodes $(0:1:-1:-1:1),(0:1:\pm i: \mp i:-1))$ and their orbits under $(\sigma, \sigma)$ are defined over
$\FF_p$.  Finally, if $p \equiv 3 (\Mod 4)$ and $p \equiv 2,3,4 (\Mod
5)$ then only the $\sigma$-nodes and $((1:1:1:1:1),(1:1:1:1:1))$ are defined over $\FF_p$.  Using the Chinese Remainder Theorem gives us the result.  
\end{proof}

\begin{prop}
If $p \equiv 1 (\Mod 20)$, all 60 nodes have blowups whose rulings are
defined over $\FF_p$.  If $p \equiv 11 (\Mod 20)$, then the
$\sigma$-nodes and $\tau$-nodes have rulings defined over $\FF_p$.  If $p \equiv 9
(\Mod 20)$, then the $\sigma$-nodes, $((1:1:1:1:1),(1:1:1:1:1))$ and
the 10 regular nodes have rulings defined over $\FF_p$.  If $p \equiv
13,17 (\Mod 20)$, then the $\sigma$-nodes and the 10 regular nodes
have rulings defined over $\FF_p$.  If $p \equiv
19 (\Mod 20)$, then the $\sigma$-nodes and $((1:1:1:1:1),(1:1:1:1:1))$
have rulings defined over $\FF_p$.  Otherwise only the $\sigma$-nodes
have rulings defined over $\FF_p$.
\end{prop}

\begin{proof} Given a node defined over $\FF_p$, the rulings
of the exceptional divisor are defined
over $\FF_p$ if and only if the determinant of the symmetric matrix
defining the node is a square in $\FF_p$.  We calculated these
determinants to be $1$, $5$ and $-12 + 16i$ (up to square factors).
Now 1 is always a square, and $-12+16i$ is a square in $\FF_p$ if $i$
is defined over $\FF_p$.  So the $\sigma$-nodes and 50 regular nodes
have rulings defined over $\FF_p$ as long as the nodes themselves are
defined over $\FF_p$.  Now 5 is a square in $\FF_p$ if and only if $p
\equiv \pm 1 (\Mod 5)$.  Using the Chinese Remainder Theorem again
gives us our result.  
\end{proof}

\subsection{Topology of $\XT$}

We need to compute the topological invariants of $\XT$.
First let $X'$ be a smooth deformation of $X$; that is, let
$X'$ be a smooth threefold obtained by intersecting $\PP^4 \times
\PP^4$ with five divisors of type $(1,1)$.  Topologically,
$\XT$ is obtained from $X'$ by contracting 60 copies
of $S^3$ and replacing them with 60 copies of $\PP^1 \times \PP^1$.
Therefore, we have 
$$\chi((\tilde X)) = \chi(X') + 60
\chi(\PP^1 \times \PP^1) = \chi(X') + 240.$$

For $0 \leq i \leq 5$, let $X^i$ denote the intersection of
$\PP^4 \times \PP^4$ with $i$ generic divisors of type $(1,1)$.  So we
have $\PP^4 \times \PP^4 = X^0$ and $X' = X^5$.  We
have the exact sequences

\begin{equation}
\begin{aligned}
0 &\longrightarrow &T(X^1) &\longrightarrow &T(\PP^4 \times \PP^4)|_{X^1}
&\longrightarrow  &N|_{X^1 / \PP^4 \times \PP^4} &\longrightarrow 0 \\
0 &\longrightarrow &T(X^i) &\longrightarrow &T(X^{i-1})|_{X^i}
&\longrightarrow &N|_{X^i / X^{i-1}} &\longrightarrow 0 \\
0 &\longrightarrow &T{X'} &\longrightarrow &T(X^4)|_{X'} &\longrightarrow
&N|_{X' / X^4} &\longrightarrow 0
\end{aligned}
\end{equation}

Let $X$ and $Y$ denote the hyperplane classes on the two copies of
$\PP^4$.  Then the Chern class of the third term in each sequence is
equal to $(1 - X - Y)$.  We thus have

\begin{equation}
c(X') = c(\PP^4 \times \PP^4) \cdot (1 - X - Y)^5.
\end{equation} 

Taking the terms of order three, we have

\begin{equation}
\begin{split}
c_3(X') &= (10 X^3 + 50 X^2Y + 50XY^2 + 10Y^3) \\
&- 5(X+Y)(10X^2+25XY+10Y^2) \\
&+ 10(X+Y)^2(5X+5Y) \\
&- 10(X+Y)^3. \\
\end{split}
\end{equation}

The Euler characteristic of $X'$ is

\begin{equation}
\chi(X') = \int_{\PP^4 \times \PP^4} c_5(N|_{{X'} / \PP^4 \times
\PP^4}) \cdot c_3(X')
\end{equation}

Now only the terms of degree 8 in the characteristic class contribute
to the integral, and the orientation class of $\PP^4 \times \PP^4$ is
$X^4 Y^4$.  The integral is thus equal to the coefficient of the $X^4 Y^4$ term in
$c_3(X') (1 + X +Y)^5$, which we compute to be $-100$.  We thus have $\chi(X') = -100$ and
$\chi(\tilde{X}) = 140$.

We can also compute most of the Hodge numbers of $\tilde{X}$.
Since $\tilde{X}$ is obtained from $X'$ by the surgery
procedure explained above, we have $h^{0,0} = h^{3,3} = 1, h^{1,0} =
h^{0,1} = h^{2,0} = h^{0,2} = 0$, and $h^{3,0} = h^{0,3} = 1$.  The
only unknown Hodge numbers are $h^{1,1} = h^{2,2}$ and $h^{1,2} =
h^{2,1}$, and we know that $2h^{1,1} - h^{2,1} = 140.$

\begin{rem}
  In section 1, we remarked that Lange
\cite{bib:Lange2} had
found a rank 2 vector bundle over $\PP^1 \times \PP^3$ whose zero
sections yielded abelian surfaces.  The Chern class of this bundle is
$1 + (2 H_1 + 4 H_3) + (8 H_1 H_3 + 6 H_3^2)$, where $H_1$ and $H_3$
denote the pullbacks to $\PP^1 \times \PP^3$ of the hyperplane class
on $\PP^1$ and $\PP^3$ respectively.  Since $c_1(\PP^1 \times \PP^3) =
(2 H_1 + 4 H_3)$, a pencil of such abelian surfaces would yield a
Calabi-Yau threefold.  A generic pencil of such abelian
surfaces would have a base locus consisting of $2 \cdot 8 \cdot
6 = 96$ nodes.  One computes the Euler characteristic of the resulting
threefold to be 176.  Since $176 < 4 \cdot 96$, it is within the realm
of possibility that resolving the singularities would yield a
threefold with $h^{1,2}$ small or zero.
\end{rem}


\section{Proof that $\XT$ is modular}

\subsection{Modular forms}

This brief review of modular forms is taken from \cite{bib:DDT}.  Recall that the special linear group $SL_2(\ZZ)$ acts on points $z$ in
the upper half plane $\mathcal{H}$:  if $\gamma = \begin{pmatrix} a &
  b \\ c & d \end{pmatrix}$, $\gamma z = \frac{a\tau+b}{c\tau+d}$.

Let $\Gamma(N)$ denote the subgroup of matrices $\gamma$ in
$SL_2(\ZZ)$ such that

\[
\gamma \equiv \begin{pmatrix} 1 & 0 \\ 0 & 1 \end{pmatrix} (\Mod N).
\]

Call a subgroup $\Gamma$ of $SL_2(\ZZ)$ a {\em congruence subgroup} if
it contains $\Gamma(N)$ for some $N$.  The {\em level} of $N$ is the smallest $N$
such that $\Gamma$ contains $\Gamma(N)$.  The most important
congruence subgroups are

\[
\Gamma_1(N) = \Biggl\{ \begin{pmatrix} a & b \\ c & d \end{pmatrix} :
\begin{pmatrix} a & b \\ c & d \end{pmatrix} \cong \begin{pmatrix} 1 &
  * \\ 0 & 1 \end{pmatrix} (\Mod N) \Biggr\},
\]

\[
\Gamma_0(N) = \Biggl\{ \begin{pmatrix} a & b \\ c & d \end{pmatrix} :
\begin{pmatrix} a & b \\ c & d \end{pmatrix} \cong \begin{pmatrix} * &
  * \\ 0 & * \end{pmatrix} (\Mod N) \Biggr\}.
\]

\begin{defn}
A {\em modular function} $f$ of weight $2k$ and level $N$ is a holomorphic
function on the upper half plane $\mathcal{H}$ such that $f(\gamma(z))
= (cz + d)^{-2k} f(z)$ for all $\gamma$ in some congruence subgroup
$\Gamma$ of level $N$.
\end{defn}

\begin{defn}
A {\em modular form} $f$ is a modular function satisfying the
following property:  for all $\gamma$ in $SL_2(\ZZ)$, the function
$(c\tau + d)^{-2k} f(\gamma\tau)$ has a Puiseux series expansion
$\sum_{n \geq 0} a_n q^{\frac{n}{h}}$ in fractional powers of $q = e^{2\pi i
  \tau}$.  We call this series the Fourier expansion of $f$ at the
cusp $\gamma^{-1}(i\infty)$, since the limit $q \rightarrow 0$
corresponds to the limit $z \rightarrow i\infty$.  A {\em cusp form}
is a modular form such that the Fourier expansion at each cusp has vanishing constant term.
\end{defn}

We will focus on congruence subgroups $\Gamma$ lying between
$\Gamma_1(N)$ and $\Gamma_0(N)$.  In this case, the matrix
$\begin{pmatrix} 1 & 1 \\ 0 & 1 \end{pmatrix}$ is in $\Gamma$.  If $f$
is a modular form with respect to $\Gamma$, we must then have
$f(\gamma z) = f(z)$ for all $z$ in $\mathcal{H}$.  Since $f$ is
periodic with period 1, its expansion at the cusp at infinity can be
written as a power series in $q = e^{2 \pi i \tau}$.

\subsection{Galois theory}

This quick review of Galois theory is taken from \cite{bib:Taylor2}.  Recall the $p$-adic absolute value on $\QQ$ given by $|x|_p = p^{-k}$
if $x$ can be expressed as $p^k \frac{a}{b}$ with $a$ and $b$ both
coprime to $p$.  The completion of the field $\QQ$ under this absolute
value gives us the field of $p$-adic numbers $\QQ_p$, and $| \quad
|_p$ extends uniquely to $\QQ_p$.  Under this metric, $\QQ_p$ is a
locally compact, totally disconnected field.  Let $G_{\QQ_p}$ denote
the absolute Galois group $\Gal(\Qpbar / \QQ_p)$; it is identified
with the group of continuous automorphisms of $\QQ_p$.  $G_{\QQ_p}$ is
a topological group; a basis at the identity is given by the
collection of subgroups of finite index.

Given an embedding of $\Qbar$ into $\Qpbar$, we obtain a closed
embedding of $G_{\QQ_p}$ into $G_{\QQ}$; this embedding varies by
conjugation as the embedding varies.

Let $\OO_{\Qpbar}$ denote the ring of integers of $\Qpbar$; it is the
ring of elements with absolute value less than or equal to 1.  It is a
local ring with maximal ideal $\mathfrak{m}_{\Qpbar}$, which is the
ideal of elements with absolute value strictly less than 1.  The
residue field $\OO_{\Qpbar} / \mathfrak{m}_{\Qpbar}$ is an algebraic
closure of $\FF_p$, which we denote by $\Fpbar$.  We thus obtain a
continuous map $G_{\QQ_p} \rightarrow G_{\FF_p}$ which is surjective.
Its kernel $I_{\QQ_p}$ is called the inertia subgroup of
$G_{\QQ_p}$.  The group $G_{\FF_p}$ is procyclic, being the inverse
limit $\lim_{\leftarrow} G(\FF_{p^k} / \FF_p)$.  The Frobenius element $Frob_p$ defined by

\begin{equation*}
\Frob_p(x) = x^p
\end{equation*}

\noindent generates a dense subgroup of $G_{\FF_p}$.

We will be looking at representations $\rho: G_{\QQ} \rightarrow
GL_d(K)$ coming from algebraic geometry in the next section.  Given an
embedding $\Qbar \hookrightarrow \Qpbar$, we have $I_{\QQ_p} \subset
G_{\QQ_p} \subset G_{\QQ}$.  We say that $\rho$ is unramified at a
prime $p$ if it is trivial on the inertia group $I_p$.  In this case,
$\rho(\Frob_p)$ is well-defined given the choice of embedding $\Qbar
\hookrightarrow \Qpbar$, and $\trace \rho(\Frob_p)$ is well-defined
independent of the embedding.

\subsection{Modularity of an algebraic variety}

Let $X$ be an projective algebraic variety defined over $\ZZ$.  We can then consider the reduction of $X$
modulo a prime $p$, i.e. the scheme $X \times_{\ZZ} \Fpbar$.

Fixing once and for all an embedding $\Qbar \subset \CC$, the Galois
group $G_{\QQ}$ acts on $X$ and thus induces an action on the
cohomology groups

\begin{equation*}
H^i(X(\CC), \QQ_l) \cong H^i_{et}(X \times_{\ZZ} \Qbar, \QQ_l).
\end{equation*}

\noindent where the left-hand side is the cohomology of $X$ in the analytic
topology.

Given a prime $p \neq l$ and an embedding $\Qbar \subset \Qpbar$, we
have the decomposition subgroup $G_{\QQ_p} \subset G_{\QQ}$.  In
addition, we have the surjective map $G_{\QQ_p} \rightarrow
G_{\FF_p}$.  It is well-known that if $p$ is a prime of good reduction
for $X$, then the representation $\rho: G_{\QQ} \rightarrow GL(H^i(X
\times_{\ZZ} \Qbar, \QQ_l))$ is unramified at $p$.  

If $p$ is a prime
of good reduction for $X$, then we can assign a well-defined value to
$\rho(\Frob_p)$ by choosing a lifting of $\Frob_p$ from $G_{\FF_p}$ to
$G_{\QQ_p}$.  By abuse of notation, such a lifting will also be
denoted by $\Frob_p$.  Furthermore, by passing from $X \times_{\ZZ} \Qpbar$
to the reduction $X \times_{\ZZ} \Fpbar$, we can consider the action of $\Frob_p$ on the cohomology
of $X \times_{\ZZ} \Fpbar$.

We say that an $n$-dimensional
algebraic variety $X$ is {\em modular} if for some subquotient $V$ of
$H^i(X, \QQ_l)$, the numbers $\Trace \Frob_p(V)$ are equal to the
coefficients of a cusp form $f$ for all but finitely many primes $p$.

En route to proving Fermat's Last Theorem, Wiles and Taylor
\cite{bib:TW} proved
that for any elliptic curve $E$ with semistable reduction at 3 and 5,
the numbers $\Trace \Frob_p (H^1(X))$ are the coefficients of a
modular form; Breuil, Conrad, Diamond and Taylor \cite{bib:BCDT} later proved the
statement for all elliptic curves $E$.

\subsection{Lefschetz, Weil and counting points}

In \'etale cohomology, we have the Lefschetz theorem \cite{bib:FK}:

\begin{thm}
If $f$ is an automorphism of the variety $X$, then the number of fixed
points of $f$ is given by the following formula:
\[
Fix(f, X) = \sum^{2n}_{i = 0} (-1)^i \Trace f^{*}(H^{i}(X)).
\]
\end{thm}

In the case $X$ is defined over $\FF_p$ and $f = \Frob_p$, $Fix(f, X)$
is simply the number of points of $X$ over $\FF_p$.  We see that the
number of points of $X$ over $\FF_p$ is related to the action of
$\Frob_p$ on the \'etale cohomology of $X$.  
 
\begin{prop}
For $p$ congruent to 1 modulo 20, the semisimplification of the action of the Frobenius map $\Frob_p$ on
$H^2(\XT \times_{\ZZ} \Fpbar, \QQ_l)$ is multiplication by $p$.
\end{prop}

\begin{proof}  We prove this statement in two steps.

\begin{step1}
The Frobenius map $\Frob_p$ acts on $H^2(X
\times_{\ZZ} \Fpbar, \QQ_l)$ by
multiplication by $p$.
\end{step1}

Consider the embedding $i: \PP^4 \times \PP^4 \rightarrow \PP^{24}$
that sends $(x,z)$ to $y$ where $y_{5i+j} = x_i z_j$.  Recall that
$X$ is a section of $\PP^4 \times \PP^4$ by five divisors of type
$(1,1)$; correspondingly $i(X)$ is a section of $i(\PP^4 \times
\PP^4)$ by five hyperplanes; as in section 3 let $X^i$ denote
successive sections of $\PP^4 \times \PP^4$ by these hyperplanes, with
the exception that here $X^i$ is $X$, instead of a deformation of $X$.
Moreover, by Bertini's theorem these sections can be chosen such that $X^{i-1} -
X^i$ is smooth for all $i$; for $p = 101$ a specific choice of
hyperplane sections is given in a Macaulay2 program in the Appendix.

By the Lefschetz theorem in \'etale cohomology, 
$$i^{*}: H^2(\PP^4 \times
\PP^4, \QQ_l) \rightarrow H^2(X \times_{\ZZ} \Fpbar, \QQ_l)$$ 
is an isomorphism that preserves the
Frobenius action.  (The proof of the Lefschetz theorem depends only on
the fact that the $i(X^{i-1} - X^i)$ are smooth and affine;
see \cite{bib:FK}, p. 106.)

Since all of $H^2(\PP^4 \times \PP^4, \QQ_l)$ can be
represented by divisors defined over $\FF_{101}$, the Frobenius map
acts by multiplication by $p$ on $H^2(\PP^4 \times \PP^4, \QQ_l)$.
Thus the Frobenius map acts likewise on $H^2(X \times_{\ZZ} \Fpbar,
\QQ_l)$.  Note that Step 1 is valid for any prime $p$, not just those
congruent to 1 modulo 20.

\begin{step2}  The semisimplification of the action of the Frobenius map $\Frob_{p}$ on
$H^2(\XT \times_{\ZZ} \Fpbar, \QQ_l)$ is multiplication by $p$.
\end{step2}

Recall that $\pi: \XT \rightarrow X$ is a blowup of
60 ordinary double points.  From the Leray spectral sequence for
$\pi$, we obtain an exact sequence

\begin{equation}
0 \longrightarrow H^2(X \times_{\ZZ} \Fpbar, \QQ_l) \longrightarrow H^2(\XT
\times_{\ZZ} \Fpbar, \QQ_l) \longrightarrow
\oplus^{60}_{i = 1} H^2(Q_i, \QQ_l),
\end{equation}

\noindent where the $Q_i$ are the exceptional divisors.  For $p \equiv
1$ (mod 20), the
rulings on all the $Q_i$ are defined over $\FF_p$.  Hence the
Frobenius map acts by multiplication by $p$ on the $Q_i$.  From the
exact sequence, we see that the semisimplification of the Frobenius map acts by multiplication
by $p$ on $H^2(\XT \times_{\ZZ} \Fpbar, \QQ_l)$.
\end{proof}

We will now concentrate our attention on the prime $p = 101$.  Let us collect the information we have so far about the cohomology of
$\XT \times_{\ZZ} \overline{\FF}_{101}$:

\begin{enumerate}
\item $h^0 = h^6 = 1.$

\item $h^1 = h^5 = 0.$

\item The semisimplification of the  $\Frob_{101}$ action on $H^2$ is
  multiplication by 101.  By Poincar\'e duality, the semisimplification
  of the $\Frob_{101}$ action on $H^4$ is multiplication by ${101}^2$.

\item $2 h^2 - h^3 = 138.$

\item $ \# X(\FF_{101}) = 1 + 101 h^2 + {101}^2 h^2 + {101}^3 - \Trace
\Frob_{101}(H^3(\XT \times_{\ZZ} \overline{\FF}_{101}, \QQ_l))$.

\item For primes $p$ of good reduction, $| \Trace \Frob_{p}(H^3(\XT
  \times_{\ZZ} \Fpbar, \QQ_l)) | \leq h^3
{p}^{\frac{3}{2}} = (138 - 2 h^2) {p}^{\frac{3}{2}}$ by the Weil conjectures.
\end{enumerate}

The number of points in $\XT$ is easily computed by the
following procedure:

\begin{enumerate}

\item For a given prime $p$, count the number of points in $G$ using a
computer.  (The code is in the Appendix.)  We start by counting points
in $G$ instead of $X$ because it is faster, the running time of
the program being an exponential function of the number of variables.

\item Add $p$ times the number of points in $E$, since each point in
$E$ is replaced by a copy of $\PP^1$ upon passage to $X$.  Note
that $E$ has points defined over $\FF_p$ only when $p \equiv 1 (\Mod
4)$.

\item Add the number of points arising from the blowup of the nodes.

\end{enumerate}

The $\sigma$-nodes are
defined for all $\FF_p$, and the rulings over the exceptional divisors
exist over all $\FF_p$.  Hence each of these five nodes adds $p^2 + p$
points to the total.

The node $((1:1:1:1:1),(1:1:1:1:1))$ is defined over all $\FF_p$, but
the other $\tau$-nodes in its orbit are defined only if $p \equiv 1 (\Mod
5)$.  The rulings over the exceptional divisors exist over $\FF_p$
only if $\sqrt{5}$ is defined in $\FF_p$, i.e. if $p \equiv \pm 1
(\Mod 5)$.  If the rulings are defined over $\FF_p$, we add $p^2 + p$
points.  Otherwise we add only $p^2$ points.

The node $((0:1:-1:-1:1),(0:1:i:-i:-1))$ is defined over $\FF_p$ if $p \equiv
1 (\Mod 4)$, as are the other nodes in its $\sigma$-orbit.  However,
the other nodes in its $H_5$-orbit are defined over $\FF_p$
only if $i$ and $\epsilon$ are both defined.  If the nodes are defined,
then the rulings over the exceptional divisors are as well.  Hence we
add $p^2 + p$ times the number of these nodes.

For $p = 101$, we obtain $\# X(\FF_{101}) = 1770940$.  Thus 

\begin{equation}
\Trace \Frob_{101}(H^3(\XT)) = 1 + 101 h^2 + {101}^2 h^2 +
{101}^3 - 1770940,
\end{equation}

\noindent so we must have

\begin{equation}
| 1 + 101 h^2 + {101}^2 h^2 + {101}^3 - 1770940 | \leq (138 - 2 h^2)
  {101}^{\frac{3}{2}}.
\end{equation}

Separating the absolute value inequality into two inequalities, we
have

\begin{equation}
\begin{aligned}
h^2 (101 + {101}^2 + 2 \cdot {101}^{\frac{3}{2}}) &\leq 138 \cdot
{101}^{\frac{3}{2}} - 1 - {101}^3 + 1770940, \\
h^2 (101 + {101}^2 - 2 \cdot {101}^{\frac{3}{2}}) &\geq -138 \cdot
{101}^{\frac{3}{2}} - 1 - {101}^3 + 1770940.
\end{aligned}
\end{equation}

Since $h^2$ must be an integer, the two inequalities force $h^2$ to be
equal to 72.  We then get $h^3 = 6$.  (This trick for computing $h^2$ is due to Werner and van
Geemen in \cite{bib:WvG}.)

For $p \not\equiv 1 (\Mod 20)$, we no longer know that the
semisimplification of $\Frob_p$ acts by
multiplication by $p$ on $H^2$.  However, we know that $\Frob_p$ acts on
$H^2(X)$ by multiplication by $p$.  We also know that
$\oplus^{60}_{i = 1} H^2(Q_i, \QQ_l)$ is spanned by algebraic cycles,
so the eigenvalues of $\Frob_p$ acting on this space are all $p$ times
roots of unity.  Using the exact sequence (29) again, the
eigenvalues of the semisimplification of $\Frob_p$ acting on
$H^2(\XT)$ are all $p$ times roots of unity.  By the Weil conjectures,
the trace of $\Frob_p$ is a rational integer.  Hence the trace of
$\Frob_p$ must be $p$ times an integer $h$.  

Suppose the eigenvalues of $\Frob_p$ acting on $H^2(\XT)$ are $p
\zeta_i$, with the $\zeta_i$ being roots of unity.  Choosing a basis
of $H^2(\XT)$ and a Poincar\'e dual basis of $H^4(\XT)$, the action of
$\Frob_p$ on $H^2(\XT)$ can be
represented as a matrix.  This matrix is similar to a matrix of upper
Jordan blocks having diagonal entries $p \zeta_i$.  By Poincar\'e
duality, the action of $\Frob_p$ acting on $H^4(\XT)$ will be $p^2$
times the contragredient of the action on $H^2(\XT)$.  Thus the matrix
of $\Frob_p$ acting on $H^4(\XT)$ will be similar to a matrix
of lower triangular blocks having diagonal entries $p^2 \overline{\zeta_i}$ with the
same multiplicities as in $H^2(\XT)$.  Hence the trace of $\Frob_p$ acting on $H^4$ is $hp^2$.

Furthermore, we have the additional piece of information that $h^2 =
72$.  So we can use the Weil conjectures again.  We now have

\begin{equation}
| 1 + ph + p^2 h + p^3 - \# (\XT) | \leq 6 p^{\frac{3}{2}}.
\end{equation}

This is equivalent to the two inequalities

\begin{equation}
h (p + p^2) \leq \# (\XT) - 1 - p^3 - 6 p^{\frac{3}{2}}
\end{equation}

\begin{equation}
h (p + p^2) \geq 1 + p^3 + 6 p^{\frac{3}{2}} - \# (\XT).
\end{equation}

It turns out that for $p = 59, 67, 71$ and $p \geq 79$, the
inequalities determine $h$ exactly.  Our calculations are listed in
Appendix B.

\begin{conj}
For primes $p$ of good reduction, the trace of $\Frob_p$ acting on
$H^2(\tilde{X})$ is $12p$, $20p$, $24p$ or $72p$ depending on
whether $i$ and $\epsilon$ are defined over $\FF_p$.
\end{conj}

There is a unique cusp form $f$ of level 5 and weight 4, whose Fourier
coefficients $a_p$ can be found in William Stein's Modular Forms
Database \cite{bib:Stein}.  It is equal to $(\eta(q) \eta(q^5))^4$, where $\eta$ is the Dedekind function.  The first few coefficients of the
$q$-expansion of $f$ are as follows:

\begin{equation}
f = q - 4 q^2 + 2 q^3 + 8 q^4 - 5 q^5 - 8 q^6 + 6 q^7 + 
\dots
\end{equation}

Comparing the values of $\Trace \Frob_p(H^3)$, we notice that for
primes $p \equiv 3 (\Mod 4)$, $a_p = \Trace \Frob_p(H^3)$, whereas for
primes $p \equiv 1 (\Mod 4)$, $\Trace \Frob_p(H^3) - a_p$ = $p (2p + 2
- \#(E))$.

\subsection{Proof of the main theorem}

We can finally prove the main theorem of this paper:

\begin{thm}
The third cohomology group $H^3$ of $\tilde{X}$ is modular in
the following sense:  as a Galois representation, its
semisimplification is a direct sum of a modular rank 2 motive $V$
and a 4-dimensional piece $W \cong H^1(\tilde{S}, \QQ_l)(-1)$.
\end{thm}

\begin{proof}  First we note that $\ST$ and $\XT$ are
defined over $\ZZ$.  Since we will only consider the coefficient group
$\QQ_l$, we can use the proper-smooth base change theorem (see for
example \cite{bib:Milne}) to pass from $\XT$
to $\XT \times_{\ZZ} \Fpbar$ (for $p \neq 2, 5, l$) and to $\XT \times_{\ZZ} \CC$.
In addition, we can pass from \'etale cohomology on $\XT \times_{\ZZ} \CC$ to
analytic cohomology by the comparison theorem.

We present the proof in several steps.

\begin{step1}  $H^1(\tilde{S} \times_{\ZZ} \Fpbar, \QQ_l)(-1) \cong \Ind_{G_{\QQ(i)}}^{G_{\QQ}}
H^1((E_1 \times_{\ZZ} \Qpbar) \times_{\Qpbar} \Fpbar, \QQ_l)(-1)$.
\end{step1}

For now, we use the analytic topology.  Recall that the map $\pi: \tilde{S}
\longrightarrow S$ is just the blowup at 50 points.  A standard Mayer-Vietoris 
argument (see for example \cite{bib:GH}, p. 473) shows that 
$H^1(\ST \times_{\ZZ} \Fpbar, \QQ_l) \cong H^1(S \times_{\ZZ} \Fpbar, \QQ_l)$.

The map $\pi: S \longrightarrow E$ is the projectivization of the rank 2
bundle $\ker L \longrightarrow E$.  Using the analytic topology for
now, the Leray-Hirsch theorem tells us that
$H^{*}(S, \QQ_l)$ is a truncated polynomial ring over $H^{*}(E, \QQ_l)$ generated by
the single element $c_1(\ker L)$, which has dimension 2.  Hence
$\pi^{*}: H^1(E, \QQ_l) \longrightarrow H^1(S, \QQ_l)$ is an isomorphism.
This statement also holds in \'etale cohomology.

Finally, recall that $E = E_1 \cup E_2$, where $E_1$ and $E_2$ are
complex conjugates of each other.  Thus $H^1(E \times_{\ZZ} \Fpbar, \QQ_l) \cong
\Ind_{G_{\QQ(i)}}^{G_{\QQ}} H^1((E_1 \times_{\ZZ} \Qpbar)
\times_{\Qpbar} \Fpbar, \QQ_l)$.  Putting everything together completes
Step 1.

\begin{step2}  $H^1(\ST \times_{\ZZ} \Fpbar, \QQ_l)(-1)$ is a subrepresentation
of $H^3(\XT \times_{\ZZ} \Fpbar, \QQ_l)$.
\end{step2}

Note that $H^1(\ST \times_{\ZZ} \Fpbar, \QQ_l)(-1)$ is isomorphic to
$H^3_{\ST}(\XT \times_{\ZZ} \Fpbar,
\QQ_l)$ (see \cite{bib:Milne}, p. 98).  So it is sufficient to show that inclusion induces an
injection $$j^{*}: H^3_{\ST}(\XT \times_{\ZZ} \Fpbar, \QQ_l) \longrightarrow
H^3(\XT \times_{\ZZ} \Fpbar, \QQ_l).$$

By base change and the comparison theorem, it suffices to prove the
same statement in the complex analytic topology.

In the analytic topology, we have the Lefschetz duality diagram (see
\cite{bib:Munkres}, p. 429)

\begin{equation*}
\begin{CD}
H^3_{\ST}(\XT \times_{\ZZ} \CC, \QQ_l) @>{j^{*}}>> H^3(\XT \times_{\ZZ} \CC, \QQ_l) \\
@VV{\cong}V                                @VV{\cong}V \\
H_3(\ST \times_{\ZZ} \CC, \QQ_l)  @>{j_*}>> H_3(\XT \times_{\ZZ} \CC, \QQ_l)
\end{CD}
\end{equation*}

\noindent where the vertical arrows are isomorphisms.  Hence it is sufficient to show that inclusion induces an injective map
$j_{*}: H_3(\ST \times_{\ZZ} \CC, \QQ_l) \longrightarrow H_3(\XT \times_{\ZZ}
\CC, \QQ_l)$.  We do this by computing
intersection classes of cycles in $H_3(\ST \times_{\ZZ} \CC, \QQ_l)$.

The 3-cycles in $\tilde{S}$ are of the form $\alpha_i \times \PP^1$ and
$\beta_j \times \PP^1$, where $\alpha_i$ and $\beta_j$ generate
$H_1(E_i)$.  Note that the $\alpha_i$ and $\beta_j$ can be chosen to
miss the exceptional divisors.

We have the fundamental result

\[
j_*(\alpha) \cap j_*(\beta) = j_*( PD[\tilde{S}]|_{\tilde{S}} \cap \alpha
\cap \beta)
\]

\noindent for any homology cycles $\alpha$ and $\beta$.

Note that

\[
[\ST] |_{\ST} = (K_{\XT} - K_{\XT} + \ST) | {\ST}
= K_{\ST} - K_{\XT} | {\ST}.
\]

But $X$ has trivial canonical bundle.  Hence $K_{\XT}$ is supported
on its exceptional fibers.  Therefore the restriction of $K_{\XT}$ to
$\ST$ is supported on the exceptional fibers of $\ST$.  Let $\gamma_i$
and $\delta_i$ be 1-cycles generating the first homology of $E_i$, and
put $\alpha_i = \gamma_i \times C_i$ and $\beta_i = \delta_i \times C_i$,
where $C_i$ is a ruling of $S_i$.

Since the cycles $\alpha_i$ and $\beta_j$ can be chosen to miss the
exceptional fibers, we have

\[
PD[\tilde{S}]|_{\tilde{S}} \cap \alpha_i \cap \beta_j = PD
[K_{\tilde{S}}] \cap \alpha_i \cap \beta_j.
\]

The canonical bundle of $\tilde{S}$ is well-known; it is simply
$-2(D) + \Sigma E_i$, where $D$ is a horizontal section and the $E_i$
are the exceptional divisors.  We also have

\[
\gamma_i \cap \gamma_j = \delta_i \cap \delta j = 0,
\]

\[
\gamma_1 \cap \delta_2 = \gamma_2 \cap \delta_1 = 0,
\]

\[
\gamma_i \cap \delta_i = C_i,
\]

\noindent where $C_i$ is a line belonging to the ruling of $S_i$. 

Therefore we have

\[
j_*(\alpha_i) \cap j_*(\alpha_j) = j_*(\beta_i) \cap j_*(\beta_j) = 0
\]

\[
j_*(\alpha_i) \cap j_*(\beta_j) = -2 \delta_{ij}.
\]

Putting these results together, we see that the intersection matrix of
the 3-cycles is 

\[
\begin{pmatrix} & -2 & & \\ 2 & & & \\ & & & -2 \\ & & 2 &
\end{pmatrix},
\]

\noindent which is nonsingular.  Hence $H_3(\tilde{S})$ injects into $H_3(\XT)$.

Upon passing to the semisimplification, we now see that as Galois representations,

\[
H^3(\XT \times_{\ZZ} \Fpbar, \QQ_l) = V \oplus \Ind_{G_{\QQ(i)}}^{G_{\QQ}} H^1(E \times_{\ZZ}
\Fpbar, \QQ_l)(-1),
\]

\noindent with $V$ some undetermined 2-dimensional piece.  This observation was
first made by Taylor \cite{bib:Taylor}.

\begin{step3}  Away from the primes of bad reduction, the
traces of $V$ coincide with the coefficients of the unique modular
form $f$ of level 5 and weight 4.
\end{step3}

For Step 3, we invoke a result of
Faltings-Serre-Livn\'e \cite{bib:Livne} which says essentially that it is enough to check the
equality at a suitably chosen finite set of primes:

\begin{thm}
(Faltings-Serre-Livn\'e) Suppose that $\rho_1$ and $\rho_2$ are two 2-adic
2-dimensional Galois representations unramified outside a set of
primes $S$.  Let $K_S$ be the smallest field containing all quadratic
extensions of $\QQ$ ramified at primes in $S$, and let $T$ be a set of
primes disjoint from $S$.  Then if

$\Trace \rho_1 \equiv \Trace \rho_2 \equiv 0$ and $\det \rho_1 \equiv
\det \rho_2 (\Mod 2)$;

$\{ \Frob_p |_{K_S} : p \in T \}$ is equal to the set $\Gal(K_S/K) -
Id$;

for all $p \in T$, $\Trace \rho_1 \Frob_p = \Trace \rho_2 \Frob_p$ and
$\det \rho_1 \Frob_p = \det \rho_2 \Frob_p$;

then $\rho_1$ and $\rho_2$ have isomorphic semisimplifications.
\end{thm}

We will apply this result in the case $\rho_1 = V$ and $\rho_2$ is the
modular Galois representation coming from the cusp form $f$ of level 5 and
weight 4.  Thus we must do the following:

Check that $\Trace \Frob_p(H^3(\XT \times_{\ZZ} \Fpbar, \QQ_l))$ is even for primes $p$ of good
reduction.

Check that the coefficients $a_p$ of the modular form $f$ are even for
$p \neq 2, 5$.

Observe that the determinants of both representations are given by
$\chi^3$, where $\chi$ is the cyclotomic character.

Find a suitable set of primes $T_s$.

Compute $\Trace \Frob_p H^3(\XT \times_{\ZZ} \Fpbar, \QQ_l)$ for all $p$ in $T_s$ and check that
these are equal to the corresponding coefficients in the modular form $f$.

\begin{lem}
$\Trace \Frob_p H^3(\XT \times_{\ZZ} \Fpbar, \QQ_l)$ is even for $p \neq 2, 5.$
\end{lem}

\begin{proof} Recall that

\begin{equation}
\Trace \Frob_p(H^3(\XT)) = p^3 + p(p+1)h + 1 - \# \XT(\FF_p).
\end{equation}

Thus we need only check that $ \# \XT(\FF_p)$ is even.

Recall the automorphism of order 4

\[
\mu = \begin{pmatrix} 1 &  &  &  & \\ & & 1 & &  \\ & & &  & 1 \\
& 1 &  & & \\ &  & & 1 & \end{pmatrix}
\]

\noindent acting on $F$, and notice that it lifts to an automorphism $(\mu,
\mu)$ of $X$.  Each point $(x,z)$ of $X$ lies in an orbit of
size 1, 2 or 4.  The points in orbits of size 2 or 4 obviously
contribute an even number of points to the total.  If any of these
points are nodes, then so are the other points in their orbits.
Blowing up each node adds either $p^2$ or $p^2 + 2p$ points to the
total.  Since there must be an even number of such nodes, blowing up
adds another even number of points to the total.  So we are left with
counting the number of points in $\XT$ coming from fixed points of
$X$.

The fixed points of $\mu$ in $F$ are of the form
$(1:x:x:x:x)$, $(0:1:1:1:1)$, $(0:1:-1:-1:1)$, and $(0:1:\pm i: \mp i:
-1)$; this is a total of $p+3$ points, which is even.  Of these
points, the only ones that are singular points in $F$ are $(1:1:1:1:1)$, $(1:0:0:0:0)$ and $(1: \pm \frac{1}{2i}:\pm \frac{1}{2i}: \pm \frac{1}{2i}: \pm \frac{1}{2i})$.  The other $p-1$ points are smooth points in $F$ and
thus lift uniquely to points in $X$ fixed by $(\mu, \mu)$.

The points $(1: \pm \frac{1}{2i}: \pm \frac{1}{2i}: \frac{1}{2i}:
\frac{1}{2i})$ (if they exist over $\FF_p$) are singular; their
inverse images under $\pi_1$ are lines.  In these lines, the only
fixed points of $(\mu, \mu)$ are

$((1:\frac{1}{2i}:\frac{1}{2i}:\frac{1}{2i}:\frac{1}{2i}),(1:\frac{1}{2i}:\frac{1}{2i}:\frac{1}{2i}:\frac{1}{2i}))$,

$((1:\frac{1}{2i}:\frac{1}{2i}:\frac{1}{2i}:\frac{1}{2i}),(0:1:-i:i:-1))$,

$((1:-\frac{1}{2i}:-\frac{1}{2i}:-\frac{1}{2i}:-\frac{1}{2i}),(1:-\frac{1}{2i}:-\frac{1}{2i}:-\frac{1}{2i}:-\frac{1}{2i}))$,

$((1:-\frac{1}{2i}:-\frac{1}{2i}:-\frac{1}{2i}:-\frac{1}{2i}),(0:1:i:-i:-1))$,

\noindent and none of these points are nodes.  So we have added another even
number of points to the total.

Finally, we look at the contribution coming from the points
$(1:1:1:1:1)$ and $(1:0:0:0:0)$ in $F$.  The fixed points of
$(\mu, \mu)$ in $X$ lying over these points are
$((1:1:1:1:1),(1:1:1:1:1))$ and $((1:0:0:0:0),(1:0:0:0:0))$.  These points
are nodes; blowing them up adds either $p^2$ or $p^2 + 2p$ points for
each node.  So
again we end up adding an even number of points to the total.  Hence
the number of points in $\XT$ is even.
\end{proof}

\begin{lem}
The coefficients $a_p$ of $f$ are all even for $p \neq 2, 5$.
\end{lem}

\begin{proof}
This proof is essentially Proposition 4.10 in \cite{bib:Livne}.  $f$ is a
modular form of level 5.  The corresponding Galois representation $\rho: \Gal(\overline{\QQ} / \QQ)
\rightarrow GL_2(\QQ_2)$ is unramified away from 2 and 5.  Let $L$ be the
extension of $\QQ$ cut out by $\ker \overline{\rho}$, where $\overline{\rho}$
denotes the reduction of $\rho$ modulo 2.  If the trace of $\rho$ were not
congruent to $0 (\Mod 2)$, then the image of $\overline{\rho}: \Gal(\overline{\QQ} / \QQ) \rightarrow GL_2(\ZZ /2)$
would contain one of the elements $\begin{pmatrix} 1 & 1 \\ 1 &
\end{pmatrix}$ or  $\begin{pmatrix} & 1 \\ 1 & 1 \end{pmatrix}$.
These matrices are of order 3.  Therefore, if $L$ is the extension of $\QQ$
cut out by $\ker \overline{\rho}$, $L$ is a $C_3$ or $S_3$ extension
of $\QQ$ unramified away from 2 and 5.  The only such field is the splitting
field of $g(x) = x^3 - x^2 + 2x + 2$ (see for instance \cite{bib:Jones}), which has Galois group $S_3$.  Since $g$
is irreducible modulo 3, the Frobenius element $\Frob_3$ has order 3 in $S_3 =
 \Gal(\overline{\QQ} / \QQ) \cong GL_2(\ZZ/2)$.  Hence the trace of
$\rho(\Frob_3) = 1$ in $\ZZ/2$.  But the Fourier coefficient $a_3$ of $f$ is
2, a contradiction.
\end{proof}

\begin{lem}
Let $T$ be the set $\{ 67, 71, 101, 103, 113, 131, 157 \}$.  Then as $p$
runs over this set, $\Frob_p$ runs over all the non-identity members
of $\Gal(\QQ_S / \QQ)$.
\end{lem}

\begin{proof}  The field $\QQ_S$ is the compositum of all
quadratic extensions of $\QQ$ unramified outside $\{2, 5 \}$, which is
$\QQ(i, \sqrt{2}, \sqrt{5})$.  Now $\Gal(\QQ_S / \QQ)$ is isomorphic
to $(\ZZ/2)^3$, and the coordinates of $\Frob_p$ are just the values of
$(\frac{-1}{p}), (\frac{2}{p})$ and $(\frac{5}{p})$.  These are
controlled by the behavior of $p$ modulo 40, as illustrated in Table 
\ref{table:legendre}.  The set of primes $T$ gives us representatives
of every element of $\Gal(\QQ_S/\QQ)$ other than the identity. 
\end{proof}

\begin{Table}
\begin{center}
\begin{tabular}{cccc}
$p (\Mod 40)$ & $(\frac{-1}{p})$ & $(\frac{2}{p})$ & $(\frac{5}{p})$ \\  
3, 27 & -1 & -1 & -1 \\
7, 23 & -1 & 1 & -1 \\
11, 19 & -1 & -1 & 1 \\
13, 37 & 1 & -1 & -1 \\  
17, 33 & 1 & 1 & -1 \\ 
21, 29 & 1 & -1 & 1 \\
31, 39 & -1 & 1 & 1 \\
1, 9 & 1 & 1 & 1 
\end{tabular}
\caption{Images of Frobenius elements in $\Gal(\QQ_s/\QQ)$.}\label{table:legendre}
\end{center}
\end{Table}

Finally, we show that for primes $p$ in the set $T$, $\Trace V$ = $\Trace \rho$.

The trace of $V$ is equal to the trace of
$\Frob_p (H^3 (\XT \times_{\ZZ} \Fpbar, \QQ_l))$ minus the trace of
$\Frob_p (H^1(\ST \times_{\ZZ} \Fpbar, \QQ_l)(-1))$.

We obtain $\Trace \Frob_p (H^3(\XT \times_{\ZZ} \Fpbar, \QQ_l))$ by counting points on $\XT$ over $\FF_p$
and by knowing what $h$ is by using the Weil conjectures, as discussed
above.  For $p$ in the set $S$, the bounds provided by the Weil
conjectures determine $h$ and thus $\Trace \Frob_p (H^3(\XT \times_{\ZZ}
\Fpbar, \QQ_l))$.  Calculations
are in the Appendix.

Recall that $H^1(\ST \times_{\ZZ} \Fpbar, \QQ_l)(-1)$ is isomorphic to $\Ind^{G_\QQ(i)}_{G_{\QQ}}
H^1(E_1 \times_{\ZZ} \Fpbar, \QQ_l)(-1)$.
For $p$ congruent to $1 (\Mod 4)$, $i$ is in $\FF_p$.  Hence $\Frob_p$ acts on $H^1(E_1 \times_{\ZZ}
\Fpbar, \QQ_l)$ and on $H^1(E_2 \times_{\ZZ} \Fpbar, \QQ_l)$,
so the the trace of $\Frob_p$ acting on $\Ind^{G_\QQ}_{G_{\QQ(i)}}
H^1(\ST \times_{\ZZ} \Fpbar, \QQ_l)(-1)$ is equal to $p$ times the sum of the traces on
$H^1(E_1 \times_{\ZZ} \Fpbar, \QQ_l)$ and $H^1(E_2 \times_{\ZZ} \Fpbar, \QQ_l)$.  In each case the trace of $\Frob_p$ acting
on $H^1(E_i \times_{\ZZ} \Fpbar, \QQ_l)$ is just $p + 1 - \# E_i (\FF_p)$, so the trace of
$\Frob_p$ acting on $H^1(\ST \times_{\ZZ} \Fpbar, \QQ_l)(-1)$ is $p (2p + 2 - \# E (\FF_p))$.

For $p$ equal to $3 (\Mod 4)$, $i$ is not in $\FF_p$.  Hence $\Frob_p$ interchanges the two
representations $H^1(E_1 \times_{\ZZ} \Fpbar, \QQ_l)$ and $H^1(E_2 \times_{\ZZ}
\Fpbar, \QQ_l)$, so the trace of
$$\Ind_{G_\QQ(i)}^{G_{\QQ}} H^1(\ST \times_{\ZZ} \Fpbar, \QQ_l)(-1)$$ is equal to zero.  

Taking the differences of these two traces gives us the traces of $V$,
which are seen to be equal to the coefficients of the modular form
$f$.
\end{proof}
\vfill
\pagebreak

\appendix


\section{C++ and Macaulay2 code}

\renewcommand{\baselinestretch}{1.62}
In this appendix we put the C++ and Macaulay2 programs we used for our computations.
\renewcommand{\baselinestretch}{0.81}

\begin{Program}
\tiny
\begin{verbatim}
#include<iostream.h> // Count points on G over Fp
int psols;

long long G(long long z0, long long z1, long long z2, long long z3,
long long z4) 
// the equation for the determinantal quintic
{
  return z0 * z0 * z0 * z1 * z4 + z1 * z1 * z1 * z2 * z0 
  + z2 * z2 * z2 * z3 * z1 + z3 * z3 * z3 * z4 * z2 + z4 * z4 * z4 * z0 * z3
  + z0 * z2 * z2 * z3 * z3 + z1 * z3 * z3 * z4 * z4 + z2 * z4 * z4 * z0 * z0
  + z3 * z0 * z0 * z1 * z1 + z4 * z1 * z1 * z2 * z2 - (z0 * z1 * z1 * z4 * z4
  + z1 * z2 * z2 * z0 * z0 + z2 * z3 * z3 * z1 * z1 + z3 * z4 * z4 * z2 * z2
  + z4 * z0 * z0 * z3 * z3 + z0 * z0 * z0 * z2 * z3 + z1 * z1 * z1 * z3 * z4
  + z2 * z2 * z2 * z4 * z0 + z3 * z3 * z3 * z0 * z1 + z4 * z4 * z4 * z1 * z2);
}

void counter(long long z0, long long z1, long long z2, long long z3, long long z4, int p) // count points on G
{
  if (G(z0, z1, z2, z3, z4) % p == 0)
    ++psols;
}

int main()
{
  long long j0, j1, j2, j3, j4;
  int p;
  cout << ``Enter a prime p: ``;
  cin >> p;
  j0 = 1; // Affine open set with j0 != 0
  for (j1 = 0; j1 < p; ++j1) {
    for (j2 = 0; j2 < p; ++j2) {
      for (j3 = 0; j3 < p; ++j3) {
        for (j4 = 0; j4 < p; ++j4) {
          counter(j0, j1, j2, j3, j4, p);
        }
      }
    }
  }
  j0 = 0;
  j1 = 1; // Affine A^3 with j0 = 0, j1 != 0 
  for (j2 = 0; j2 < p; ++j2) {
    for (j3 = 0; j3 < p; ++j3) {
      for (j4 = 0; j4 < p; ++j4) {
        counter(j0, j1, j2, j3, j4, p);
      }
    }
  }
  j1 = 0;
  j2 = 1; // Affine A^2 with j0 = j1 = 0, j2 != 0
  for (j3 = 0; j3 < p; ++j3) {
    for (j4 = 0; j4 < p; ++j4) {
     counter(j0, j1, j2, j3, j4, p);
    }
  }
  j2 = 0;
  j3 = 1; // Affine A^1 with j0 = j1 = j2 = 0, j3 != 0
  for (j4 = 0; j4 < p; ++j4) {
    counter(j0, j1, j2, j3, j4, p);
  }
  counter(j0, j1, j2, 0, 1, p);
  cout << ``\n''
    ``The number of solutions mod p is `` << psols << ``\n'';
}
\end{verbatim}
\normalsize
\caption{Counting points on $G$.}
\end{Program}

\pagebreak

\begin{Program}
\tiny
\begin{verbatim}

#include <iostream.h> // Count points on E modulo a prime p

int psols;

long long E0(long long z0, long long z1, long long z2, long long z3, long long z4, long long y)
{
  return (-y * z0 * z0 - z1 * z4 + y * y * z2 * z3);
}

long long E1(long long z0, long long z1, long long z2, long long z3, long long z4, long long y)
{
  return (-y * z1 * z1 - z2 * z0 + y * y * z3 * z4);
}

long long E2(long long z0, long long z1, long long z2, long long z3, long long z4, long long y)
{
  return (-y * z2 * z2 - z3 * z1 + y * y * z4 * z0);
}

long long E3(long long z0, long long z1, long long z2, long long z3, long long z4, long long y)
{
  return (-y * z3 * z3 - z4 * z2 + y * y * z0 * z1);

long long E4(long long z0, long long z1, long long z2, long long z3, long long z4, long long y)
{
  return (-y * z4 * z4 - z0 * z3 + y * y * z1 * z2);
}

void counter(long long z0, long long z1, long long z2, long long z3,
long long z4, long long y, int p) 
// count points
{
  if (E0(z0, z1, z2, z3, z4, y) % p == 0) {
    if (E1(z0, z1, z2, z3, z4, y) % p == 0) {
      if (E2(z0, z1, z2, z3, z4, y) % p == 0) {
        if (E3(z0, z1, z2, z3, z4, y) % p == 0) {
          if (E4(z0, z1, z2, z3, z4, y) % p == 0) {
            ++psols;
          }
        }
      }
    }
  }
}

int main()
{
  long long j0, j1, j2, j3, j4, y;
  int p;
  cout << ``Enter a parameter y: ``; // usually y is a square root of -1
  cin >> y;
  cout << ``Enter a prime p: ``;
  cin >> p;
  j0 = 1; // Affine open set with j0 != 0
  for (j1 = 0; j1 < p; ++j1) {
    for (j2 = 0; j2 < p; ++j2) {
      for (j3 = 0; j3 < p; ++j3) {
        for (j4 = 0; j4 < p; ++j4) {
          counter(j0, j1, j2, j3, j4, y, p);
        }
      }
    }
  }
  j0 = 0;
  j1 = 1; // Affine A^3 with j0 = 0, j1 != 0
  for (j2 = 0; j2 < p; ++j2) {
    for (j3 = 0; j3 < p; ++j3) {
      for (j4 = 0; j4 < p; ++j4) {
        counter(j0, j1, j2, j3, j4, y, p);
      }
    }
  }
  j1 = 0;
  j2 = 1; // Affine A^2 with j0 = j1 = 0, j2 != 0
  for (j3 = 0; j3 < p; ++j3) {
    for (j4 = 0; j4 < p; ++j4) {
      counter(j0, j1, j2, j3, j4, y, p);
    }
  }
  j2 = 0;
  j3 = 1; // Affine A^1 with j0 = j1 = j2 = 0, j3 != 0
  for (j4 = 0; j4 < p; ++j4) {
    counter(j0, j1, j2, j3, j4, y, p);
  }
  counter(j0, j1, j2, 0, 1, y, p);
  cout << ``\n''
    ``The number of solutions to Ey mod p is `` << psols << ``\n'';
\end{verbatim}
\normalsize
\caption{Counting points on $E_1$ and $E_2$.}
\end{Program}

\pagebreak

\begin{Program}
\tiny
\begin{verbatim}
F = ZZ;

R = F[x0,x1,x2,x3,x4,z0,z1,z2,z3,z4];

f0 = 0 * 0  -x3*z1 + x1*z2 + x4*z3 - x2*z4;
f1 = -x3*z0 +0 * 0 - x4*z2 + x2*z3 + x0*z4;
f2 =  x1*z0 -x4*z1 + 0 * 0 - x0*z3 + x3*z4;
f3 =  x4*z0 +x2*z1 - x0*z2 + 0 * 0 - x1*z4;
f4 = -x2*z0 +x0*z1 + x3*z2 - x1*z3 + 0 * 0;

L = matrix{{0,z2,-z4,-z1,z3},{z4,0,z3,-z0,-z2},{-z3,z0,0,z4,-z1},{-z2,-z4,z1,0,z0},{z1,-z3,-z0,z2,0}};
G = minors(4,L);

M = matrix{{0,-x3,x1,x4,-x2},{-x3,0,-x4,x2,x0},{x1,-x4,0,-x0,x3},{x4,x2,-x0,0,-x1},{-x2,x0,x3,-x1,0}};
F = minors(4,M);

S0 = diff f0;
S1 = S0 || diff f1;
S2 = S1 || diff f2;
S3 = S2 || diff f3;
S = S3 || diff f4;

X = ideal(f0,f1,f2,f3,f4);
sing = minors(5,S);
slocus = X + sing;
slocus2 = slocus + G + F;
k = 100*x0*x1*x2*x3*x4*z0*z1*z2*z3*z4;

--now check that k is contained in slocus2
\end{verbatim}
\normalsize
\caption{Proving that nodes on $X$ have $x_0 = 0$ or $z_0 = 0$.}
\end{Program}

\pagebreak

\begin{Program}
\tiny
\begin{verbatim}
-- Let G be the quintic threefold as in our thesis.  E1 U E2 is a union of two elliptic curves.
-- We show that over any field of characteristic not equal to 2,
-- the set of points in E1 U E2 is precisely the set of rank 3 points of G.

-- I is the ideal of the singular locus of G.

F = ZZ;
R = F[z0,z1,z2,z3,z4];
L = matrix{{0,-z4,z3,z2,-z1},{-z2,0,-z0,z4,z3},{z4,-z3,0,-z1,z0},{z1,z0,-z4,0,-z2},{-z3,z2,z1,-z0,0}};
I = minors(4,L);

-- the ideal of E1 is generated by the elements I ri + si, i = 0,1,2,3,4.
-- the ideal of E2 is generated by the elements -I ri + si, i = 0,1,2,3,4.

r0 = z0^2;
r1 = z1^2;
r2 = z2^2;
r3 = z3^2;
r4 = z4^2;
s0 = (z1*z4 + z2*z3);
s1 = (z2*z0 + z3*z4);
s2 = (z3*z1 + z4*z0);
s3 = (z4*z2 + z0*z1);
s4 = (z0*z3 + z1*z2);

-- if the characteristic of our base field is not 2, then the ideal of E1 * E2 is generated by the elements
-- pi = (I ri + si)*(-I ri + si) = (ri^2 + si^2)
-- and the elements (I ri + si)*(-I rj + sj) = ri*rj + si*sj + I(ri*sj - si*rj).
-- But since the characteristic is not 2, the generators ri*rj + si*sj + I(ri*sj - si*rj)
-- and ri*rj + si*sj - I(ri*sj - si*rj)
-- can be replaced by ri*rj + si*sj and ri*sj - si*rj.

p0 = r0^2 + s0^2;
p1 = r1^2 + s1^2;
p2 = r2^2 + s2^2;
p3 = r3^2 + s3^2;
p4 = r4^2 + s4^2;

q01 = r0 * r1 + s0 * s1;
q10 = -r0 * s1 + r1 * s0;
q02 = r0 * r2 + s0 * s2;
q20 = -r0 * s2 + r2 * s0;
q03 = r0 * r3 + s0 * s3;
q30 = -r0 * s3 + r3 * s0;
q04 = r0 * r4 + s0 * s4;
q40 = -r0 * s4 + r4 * s0;
q12 = r1 * r2 + s1 * s2;
q21 = -r1 * s2 + r2 * s1;
q13 = r1 * r3 + s1 * s3;
q31 = -r1 * s3 + r3 * s1;
q14 = r1 * r4 + s1 * s4;
q41 = -r1 * s4 + r4 * s1;
q23 = r2 * r3 + s2 * s3;
q32 = -r2 * s3 + r3 * s2;
q24 = r2 * r4 + s2 * s4;
q42 = -r2 * s4 + r4 * s2;
q34 = r3 * r4 + s3 * s4;
q43 = -r3 * s4 + r4 * s3;

E = ideal(p0,p1,p2,p3,p4,q01,q10,q02,q20,q03,q30,q04,q40,q12,q21,q13,q31,q14,q41,q23,q32,q24,q42,q34,q43);

-- I is the ideal of the rank 3 locus of G.
-- E is the ideal of E1 U E2.
-- we show that the radicals of these ideals are the same.

Erad = radical E;
Irad = radical I;
Egens = transpose gens gb Erad;
Igens = transpose gens gb Irad;

-- one checks that the elements of Egens and the elements of Igens are identical, up to factors of 2.
-- Hence over a base field of characteristic not equal to 2, Egens and Igens cut out the same points.
-- The singular locus of G thus consists of E and some extra rank 4 points.
\end{verbatim}
\normalsize
\caption{Proof that $E_1 \cup E_2$ is the rank 3 locus of $G$.}
\end{Program}

\pagebreak

\begin{Program}
\tiny
\begin{verbatim}

-- We want to find the singular points of G.
-- with z0 = 0

F = ZZ;
R = F[z0,z4,z2,z3,z1];

L = matrix{{0,z2,-z4,-z1,z3},{z4,0,z3,-z0,-z2},{-z3,z0,0,z4,-z1},{-z2,-z4,z1,0,z0},{z1,-z3,-z0,z2,0}};

-- here G is actually 2 times what it should be but that's ok
-- if we are not in characteristic 2

G = det L;

H = diff(G);

-- the ideal slocus defines the singular locus of G

slocus = ideal H + ideal (G);
slocusgens = transpose gens gb slocus;

-- the ideal slocusz0 defines the set of singular points of G
-- with z0 = 0

slocusz0 = slocus + ideal(z0);
slocusz0gens = transpose gens gb slocusz0;

-- checking the Grobner basis slocusz0gens, we find that either some
-- other zi is 0, in which case z is a sigma-node, or
-- z1 z4 + z2 z3 = 0.  In this case,

slocus2z0 = slocusz0 + ideal(z1 * z4 + z2 * z3);
slocus2z0gens = transpose gens gb slocus2z0;

-- checking the Grobner basis slocus2z0gens, 
-- and setting z1 = 1, one concludes that
-- z0 = 0, z1 = 1, z2 = z3^3, z4 = -z3^4, 
-- and z3^10 = -1.  
\end{verbatim}
\normalsize
\caption{Finding nodes of G that satisfy $z_0 = 0$.}
\end{Program}

\begin{Program}
\tiny
\begin{verbatim}
-- consider the segre variety P4 x P4 in P24.  our threefold X is a section of P4 by five
-- hyperplanes.  we show that one can choose five hyperplanes such that no singularities appear until 
-- the final slice.  we can then apply the Lefschetz theorem to compute Hodge numbers of X.

F = ZZ/101;
R = F[x0, x1, x2, x3, z0, z1, z2, z3];
x4 = 1;
z4 = 1;
f0 = -x3 * z1 + x1 * z2 + x4 * z3 - x2 * z4;
f1 = -x3 * z0 - x4 * z2 + x2 * z3 + x0 * z4;
f2 = x1 * z0 - x4 * z1 - x0 * z3 + x3 * z4;
f3 = x4 * z0 + x2 * z1 - x0 * z2 - x1 * z4;
f4 = -x2 * z0 + x0 * z1 + x3 * z2 - x1 * z3;
X0 = ideal(f0+3*f1+7*f2+2*f3+10*f4);
X1 = X0 + ideal(f0+5*f1+27*f2+53*f3+18*f4);
X2 = X1 + ideal(f0+54*f1+33*f2+42*f3+20*f4);
X3 = X2 + ideal(f0+9*f2+38*f2+19*f3+64*f4);
X4 = X3 + ideal(f0);
J0 = diff (f0+3*f1+7*f2+2*f3+10*f4);
J1 = J0 || diff (f0+5*f1+27*f2+53*f3+18*f4);
J2 = J1 || diff (f0+54*f1+33*f2+42*f3+20*f4);
J3 = J2 || diff (f0+9*f2+38*f2+19*f3+64*f4);
J4 = J3 || diff f0;
sing0 = minors(1,J0);
sing1 = minors(2,J1);
sing2 = minors(3,J2);
sing3 = minors(4,J3);
sing4 = minors(5,J4);
slocus0 = X0 + sing0;
slocus1 = X1 + sing1;
slocus2 = X2 + sing2;
slocus3 = X3 + sing3;
slocus4 = X4 + sing4;

-- check that dim slocus0 = dim slocus1 = dim slocus2 = dim slocus 3 = -1
-- and dim slocus4 = 0

\end{verbatim}
\normalsize
\caption{Proof that $\PP^4 \times \PP^4$ can be successively sliced to
  obtain $X$ in such a manner that singularities only appear at
  the end.}
\end{Program}

\pagebreak

\section{Calculating traces of $\Frob_p$}

\begin{Table}[h]
\scriptsize
\begin{center}
\begin{tabular}{lcccccc}
$p$ & 59 & 67 & 71 & 79 & 83 & 89 \\
\\
$\#G(\FF_p)$& 225766 & 327706 & 407910 & 529886 & 613006 & 751756 \\
$\sigma$-nodes defined over $\FF_p$ & 5 & 5 & 5 &
5 & 5 & 5 \\
$\tau$-nodes defined over $\FF_p$ & 1 & 1 & 5 &
1 & 1 & 1 \\
Other nodes defined over $\FF_p$ & 0 & 0 & 0 & 0 & 0 & 10 \\
Points on $E_1 \cup E_2$ & 0 & 0 & 0 & 0 & 0 & 180 \\
\\
$i$ in $\FF_p$? & 0 & 0 & 0 & 0 & 0 & 1 \\
$\sqrt{5}$ in $\FF_p$? & 1 & 0 & 1 & 1 & 0 & 1 \\
$\epsilon$ in $\FF_p$? & 0 & 0 & 1 & 0 & 0 & 0 \\
\\
$\# X(\FF_p)$ & 247360 & 355310 & 459740 & 568280 & 655170 &
897360\\
\\
$p^3 + 1 - \# X(\FF_p)$ & -41980 & -54546 & -101828 & -75240 &
-83382 & -192390 \\
$p^2 + p$ & 3540 & 4556 & 5112 & 6320 & 6972 & 8010 \\
$h$ & 12 & 12 & 20 & 12 & 12 & 24 \\
$\trace \Frob_p$ on $H^3$ & 500 & 126 & 412 & 600 & 282 & -150 \\
\\
$6 p^{\frac{3}{2}}$ & 2719.2 & 3290.6 & 3589.6 & 4213.1 & 4537.0 &
5037.8 \\
$a_p$ & 500 & 126 & 412 & 600 & 282 & -150
\\
$\trace \Frob_p - a_p$ & 0 & 0 & 0 & 0 & 0 & 0 \\
$(\trace \Frob_p - a_p)/p$ & 0 & 0 & 0 & 0 & 0 & 0 \\
\\
$2p + 2 - \# (E_1 \cup E_2)(\FF_p)$ & & & & & & 0 \\

\end{tabular}
\normalsize
\caption{Counting points on $X$ and calculating traces of $\Frob_p$.}\label{table:traces}
\end{center}
\end{Table}

\begin{Table}
\scriptsize
\begin{center}
\begin{tabular}{lcccccc}
$p$ & 97 & 101 & 103 & 107 & 109 & 113 \\
\\
$\#G(\FF_p)$& 967966 & 1126560 & 1157186 & 1295146 & 1365776 & 1517046 \\
$\sigma$-nodes defined over $\FF_p$ & 5 & 5 & 5 &
5 & 5 & 5 \\
$\tau$-nodes defined over $\FF_p$ & 1 & 5 & 1 &
1 & 1 & 1 \\
Other nodes defined over $\FF_p$ & 10 & 50 & 0 & 0 & 10 & 10 \\
Points on $E_1 \cup E_2$ & 170 & 200 & 0 & 0 & 220& 230 \\
\\
$i$ in $\FF_p$? & 1 & 1 & 0 & 0 & 1 & 1 \\
$\sqrt{5}$ in $\FF_p$? & 0 & 1 & 0 & 0 & 1 & 0 \\
$\epsilon$ in $\FF_p$? & 0 & 1 & 0 & 0 & 0 & 0 \\
\\
$\# X(\FF_p)$ & 1137910 & 1770940 & 1221870 & 1364910 & 1583340 &
1750730\\
\\
$p^3 + 1 - \# X(\FF_p)$ & -225236 & -740638 & -129142 & -139866 &
-288310 & -307832 \\
$p^2 + p$ & 9506 & 10302 & 10712 & 11556 & 11990 & 12882 \\
$h$ & 24 & 72 & 12 & 12 & 24 & 24 \\
$\trace \Frob_p$ on $H^3$ & 2908 & 1106 & -598 & -1194 & -550 & 1336 \\
\\
$6 p^{\frac{3}{2}}$ & 5732.1 & 6090.3 & 3589.6 & 6272.1 & 6640.9 &
6828.0 \\
$a_p$ & 386 & 702 & -598 & -1194 & -550 & 1562
\\
$\trace \Frob_p - a_p$ & 2522 & 404 & 0 & 0 & 0 & -226 \\
$(\trace \Frob_p - a_p)/p$ & 26 & 4 & 0 & 0 & 0 & -2 \\
\\
$2p + 2 - \# (E_1 \cup E_2)(\FF_p)$ & 26 & 4 & & &  0 & -2 \\
\end{tabular}
\normalsize
\caption{Counting points on $X$ and calculating traces of
  $\Frob_p$, cont.}\label{table:traces2}
\end{center}
\end{Table}

\pagebreak

\begin{Table}
\scriptsize
\begin{center}
\begin{tabular}{lcccccc}
$p$ & 127 & 131 & 137 & 139 & 149 & 151 \\
\\
$\#G(\FF_p)$& 2143566 & 24219190 & 2685206 & 2802246 & 3437616 &
3669110  \\
$\sigma$-nodes defined over $\FF_p$ & 5 & 5 & 5 &
5 & 5 & 5 \\
$\tau$-nodes defined over $\FF_p$ & 1 & 5 & 1 &
1 & 1 & 5 \\
Other nodes defined over $\FF_p$ & 0 & 0 & 10 & 0 & 10 & 0  \\
Points on $E_1 \cup E_2$ & 0 & 0 & 290 & 0 & 300 & 0  \\
\\
$i$ in $\FF_p$? & 0 & 0 & 1 & 0 & 1 & 0  \\
$\sqrt{5}$ in $\FF_p$? & 0 & 1 & 0 & 1 & 1 & 1  \\
$\epsilon$ in $\FF_p$? & 0 & 1 & 0 & 0 & 0 & 1  \\
\\
$\# X(\FF_p)$ & 2241610 & 2596140 & 3029350 & 2919840 & 3842300 &
3900140 \\
\\
$p^3 + 1 - \# X(\FF_p)$ & -193226 & -348048 & -457966 & -234220 &
-534350 & -457188 \\
$p^2 + p$ & 16256 & 17292 & 18906 & 19460 & 22350 & 22952 \\
$h$ & 12 & 20 & 24 & 12 & 24 & 20  \\
$\trace \Frob_p$ on $H^3$ & 1846 & -2208 & -4252 & -700 & 2050 & 1852
 \\
\\
$6 p^{\frac{3}{2}}$ & 8587.4 & 8996.2 & 9621.3 & 9832.8 & 10912.7 &
11133.2 \\
$a_p$ & 1846 & -2208 & -2334 & -700 & 2050 & 1852
\\
$\trace \Frob_p - a_p$ & 0 & 0 & -1918 & 0 & 0 & 0  \\
$(\trace \Frob_p - a_p)/p$ & 0 & 0 & -14 & 0 & 0 & 0  \\
\\
$2p + 2 - \# (E_1 \cup E_2)(\FF_p)$ & & & -14 & & 0 &  \\
\end{tabular}
\normalsize
\caption{Counting points on $X$ and calculating traces of
  $\Frob_p$, cont.}\label{table:traces3}
\end{center}
\end{Table}

\begin{Table}
\scriptsize
\begin{center}
\begin{tabular}{lc}
$p$ & 157 \\
\\
$\#G(\FF_p)$ & 4019026 \\
$\sigma$-nodes defined over $\FF_p$ & 5 \\
$\tau$-nodes defined over $\FF_p$  & 1 \\
Other nodes defined over $\FF_p$  & 10 \\
Points on $E_1 \cup E_2$ & 350 \\
\\
$i$ in $\FF_p$?  & 1 \\
$\sqrt{5}$ in $\FF_p$?  & 0 \\
$\epsilon$ in $\FF_p$? & 0 \\
\\
$\# X(\FF_p)$ & 4473070 \\
\\
$p^3 + 1 - \# X(\FF_p)$ & -603176 \\
$p^2 + p$  & 24806 \\
$h$ & 24 \\
$\trace \Frob_p$ on $H^3$ & -7832 \\
\\
$6 p^{\frac{3}{2}}$  & 11803.3 \\
$a_p$ & -2494
\\
$\trace \Frob_p - a_p$  & -5338 \\
$(\trace \Frob_p - a_p)/p$  & -34 \\
\\
$2p + 2 - \# (E_1 \cup E_2)(\FF_p)$ & -34 \\
\end{tabular}
\normalsize
\caption{Counting points on $X$ and calculating traces of
  $\Frob_p$, cont.}\label{table:traces4}
\end{center}
\end{Table}
\pagebreak

\end{document}